\documentclass[opre,nonblindrev]{informs3}

\OneAndAHalfSpacedXII 
\usepackage{endnotes}
\let\footnote=\endnote
\let\enotesize=\normalsize
\def\notesname{Endnotes}%
\def\makeenmark{\hbox to1.275em{\theenmark.\enskip\hss}}
\def\enoteformat{\rightskip0pt\leftskip0pt\parindent=1.275em
  \leavevmode\llap{\makeenmark}}
\usepackage[utf8]{inputenc}
\usepackage[ruled,vlined]{algorithm2e}
\usepackage{amsmath}
\usepackage{amssymb}
\usepackage{bbm}
\usepackage{glossaries}
\usepackage{graphics}
\usepackage{mathtools}
\usepackage{nomencl}
\usepackage{subcaption}

\newcommand{\stackedcell}[2][c]{%
  \begin{tabular}[#1]{@{}c@{}}#2\end{tabular}}

\newacronym{ad}{AD}{automatic differentiation}
\newacronym{cart}{CART}{classification and regression trees}
\newacronym{doe}{DoE}{design of experiments}
\newacronym{gp}{GP}{Gaussian process}
\newacronym{iai}{IAI}{Interpretable AI}
\newacronym{knn}{$k$NN}{$k$-Nearest Neighbors}
\newacronym{kkt}{KKT}{Karush-Kuhn-Tucker}
\newacronym{nn}{NN}{Nearest Neighbors}
\newacronym{lh}{LH}{Latin hypercube}
\newacronym{lp}{LP}{linear program}
\newacronym{mi}{MI}{mixed-integer}
\newacronym{mio}{MIO}{mixed-integer optimization}
\newacronym{milo}{MILO}{mixed-integer linear optimization}
\newacronym{minlp}{MINLP}{mixed-integer nonlinear program}
\newacronym{ml}{ML}{machine learning}
\newacronym{nlp}{NLP}{nonlinear program}
\newacronym{oa}{OA}{orthogonal array}
\newacronym{oct}{OCT}{optimal classification tree}
\newacronym{octs}{OCTs}{optimal classification trees}
\newacronym{octh}{OCT-H}{optimal classification tree with hyperplanes}
\newacronym{octhagon}{OCT-HaGOn}{OCT-H for Global Optimization}
\newacronym{olh}{OLH}{optimal Latin hypercube}
\newacronym{oos}{OOS}{on-orbit servicing}
\newacronym{ort}{ORT}{optimal regression tree}
\newacronym{orts}{ORTs}{optimal regression trees}
\newacronym{orth}{ORT-H}{optimal regression tree with hyperplanes}
\newacronym[]{pgd}{PGD}{projected gradient descent}
\newacronym{pwl}{PWL}{piecewise linear}
\newacronym{soc}{SOC}{second order conic}
\newacronym{wlog}{w.l.o.g.}{without loss of generality}

\newcommand{\R}{\mathbb{R}}
\newcommand{\I}{\mathbb{I}}

\newcommand{\bd}{{\bf d}}

\newcommand{\bg}{{\bf g}}
\newcommand{\bh}{{\bf h}}

\newcommand{\bu}{{\bf u}}

\newcommand{\bx}{{\bf x}}
\newcommand{\by}{{\bf y}}
\newcommand{\bz}{{\bf z}}

\newcommand{\bbp}{{\bf P}}

\newcommand{\balpha}{{\boldsymbol{\alpha}}}

\usepackage{natbib}
\bibpunct[, ]{(}{)}{,}{a}{}{,}%
\def\bibfont{\small}%
\TheoremsNumberedThrough     % Preferred (Theorem 1, Lemma 1, Theorem 2)
\ECRepeatTheorems

\EquationsNumberedThrough    % Default: (1), (2), ...
\begin{document}
%%%%%%%%%%%%%%%%

%\setcounter{page}{1}
%\DOI{10.1287/xxxx.0000.0000}%

\RUNAUTHOR{Bertsimas and \"Ozt\"urk}

\RUNTITLE{Global Optimization via Optimal Decision Trees}

% \TITLE{Bundling Information Goods of Decreasing Value}
\TITLE{Global Optimization via Optimal Decision Trees}

\ARTICLEAUTHORS{%
\AUTHOR{Dimitris Bertsimas}
\AFF{Sloan School of Management, Massachusetts Institute of Technology, Cambridge, MA, USA, \EMAIL{dbertsim@mit.edu}} %, \URL{}}
\AUTHOR{Berk \"Ozt\"urk}
\AFF{Department of Aeronautics and Astronautics, Massachusetts Institute of Technology, Cambridge, MA, USA, \EMAIL{bozturk@mit.edu}}
} % end of the block

\ABSTRACT{
The global optimization literature places large emphasis on reducing intractable optimization problems into more tractable structured optimization forms. In order to achieve this goal, many existing methods are restricted to optimization over explicit constraints and objectives that use a subset of possible mathematical primitives. These are limiting in real-world contexts where more general explicit and black box constraints appear. Leveraging the dramatic speed improvements in mixed-integer optimization (MIO) and recent research in machine learning, we propose a new method to learn MIO-compatible approximations of global optimization problems using optimal decision trees with hyperplanes (OCT-Hs). This constraint learning approach only requires a bounded variable domain, and can address both explicit and inexplicit constraints. We solve the MIO approximation efficiently to find a near-optimal, near-feasible solution to the global optimization problem. We further improve the solution using a series of projected gradient descent iterations. We test the method on a number of numerical benchmarks from the literature as well as real-world design problems, demonstrating its promise in finding global optima efficiently.
}

% \HISTORY{This paper was first submitted in November, 2021.}

\maketitle
\section{Introduction}

Global optimization seeks to address the following problem
\begin{align}
    \begin{split}
        \label{eq:globalopt}
            \underset{\bx}{\text{min}} &~~f(\bx) \\
            \text{s.t.} &~~\bg(\bx) \geq \mathbf{0}, \\
            &~~\bh(\bx) = \mathbf{0}, \\
            &~~~\bx \in \mathbb{Z}^m \times \R^{n-m},
    \end{split}
\end{align}
where $f$, $\bg$ and $\bh$ are the objective function, inequality constraints and equality constraints, respectively, and $\bx$ is a vector of decision variables. The objective functions and constraints may or may not conform to any specific mathematical structure, unlike linear or convex optimization problems, and variables can be continuous or integer. 

Existing global optimizers approximate problem \eqref{eq:globalopt} into forms compatible with efficient optimization. These optimizers use three major approaches, which are gradient-based methods, outer approximations, and \gls{mio} methods. The gradient-based approach is used by popular nonlinear solvers such as CONOPT and IPOPT. These solvers initialize their solution procedure using feasible solutions found via efficient heuristics. Then, they solve a series of gradient descent iterations, confirming optimality via satisfaction of the \gls{kkt} conditions. As detailed by \cite{ArneS.Drud1994}, CONOPT relies on a generalized reduced-gradient algorithm, linearizing the constraints and solving a sequence of linear-searching gradient steps, maintaning feasibility to tolerance at each step. \cite{Wachter2006} describe IPOPT's primal-dual barrier approach. It relaxes the constrained global optimization problem into an unconstrained optimization problem using a logarithmic barrier function, then uses a damped Newton's method to reduce the optimality gap to a desired tolerance. These gradient-based optimizers are efficient and effective in the presence of nonlinear constraints that are sparse, being able to solve problems on the order of 1000 variables and constraints on unremarkable personal computers in minutes to local optimality. 

Another approach is an outer approximation approach, described by \cite{Horst1989}. This approach simplifies a global optimization problem by approximating constraints via linear and nonlinear cuts that preserve the original the feasible set over decision variables $\bx$. This approach is effective for constraints 
with certain mathematical structure (e.g., linearity of integer variables and convexity of nonlinear functions considered by \cite{Duran1986}, or concavity or bilinearity of constraints considered by \cite{Bergamini2008}), where mathematically efficient outer approximators exist. While these approaches are effective, they have found less commercial success due to their problem specific nature. 

A final approach, and one that meshes naturally with optimization over integer variables, couples \gls{mio} with outer approximations. \cite{Ryoo1996} present an impactful approach called the branch-and-reduce method, which relies on recursively partitioning the domain of each constraint and objective over the decision variables, and bounding their values in each subdomain by examining their mathematical primitives. Such recursive partitioning creates a branch-and-bound tree, the solution to which has guarantees of global optimality through the bounding and pruning process inherent in solving \gls{mio} problems via branch-and-bound. This method has seen success in BARON [\cite{sahinidis:baron:21.1.13}], a popular commercial global optimizer. 

While the aforementioned approaches are effective in addressing certain classes of global optimization problems, each of these approaches has weaknesses. In general, gradient-based approaches rely on good initial feasible solutions, and are ineffective in presence of integer decision variables. Outer approximation approaches fail to generalize to global optimization problems with general nonlinearities. While being more general than outer approximation methods, existing \gls{mio} approaches don't scale as well due to their combinatorial nature.

Perhaps more importantly, in pursuit of mathematical efficiency, many global optimizers place additional constraints on the forms of constraints, requiring constraints to use a small subset of possible mathematical primitives. For example, BARON ``can handle functions that involve $\mathrm{exp}(x)$, $\mathrm{ln}(x)$, $x^{\alpha}$ for real $\alpha$, and $\beta^x$ for real $\beta$" [\cite{sahinidis:baron:21.1.13}]. Constraints from the real world do not always adhere to these forms, and often involve other classes of functions such as trigonometric functions, signomials, and piecewise-discontinuous functions. It is often not possible to transform these functions into forms compatible with existing global optimizers. These optimizers face even greater challenges when dealing with objectives and constraints that are black box. Black box constraints are \emph{inexplicit}, meaning that they have no analytical representations, such as when constraints are the outcomes of simulations. 

In this paper, we propose a new approach to reformulate global optimization problems as \gls{mio} problems using \gls{ml}, leveraging work by \cite{Bertsimas2017,Bertsimas2019} on the \gls{octh} and the \gls{orth}. The approach addresses global optimization with arbitrary explicit and inexplicit constraints. The only requirement for the proposed method is a bounded feasible domain for the subset of decision variables $\bx$ present in nonlinear constraints. 

In our proposed method, we approximate each constraint that is outside of the scope of efficient mathematical optimization using an \gls{octh}. More specifically, each nonlinear constraint $g_i(\bx) \geq 0$ is approximated by an \gls{octh} $T_{i}$ trained on data $\{(\tilde{\bx}_{k}, \I(g_i(\tilde{\bx}_{k}) \geq 0)),~k \in [n]\}$, where $\tilde{\bx}_k$ is an outcome of decision variables, $\I$ is the indicator function, and $g_i(\tilde{\bx}_k)$ is the left-hand-side of the constraint evaluated at $\tilde{\bx}_k$. Thus, tree $T_{i}$ makes an approximation of the feasible space of constraint $g_i(\bx) \geq 0$, predicting (with some error)  whether an outcome of decision variables satisfies the constraint. This approach also extends to approximate each nonlinear equality $h_j(\bx) = 0$, and approximates nonlinear objective functions via \gls{orth}s. 

The approximating trees allow for a natural \gls{mio} approximation of the underlying constraints. Each feasible leaf of an \gls{octh} is reached by a decision path defining an intersection of halfspaces, i.e. a polyhedron. Constraints may thus be approximated as a union of feasible polyhedra of the approximating \gls{octh}s using disjunctive constraints. We solve this efficient \gls{mio} approximation of the original problem to obtain a near-feasible and near-optimal solution, and then use gradient-based methods to repair the solution to be feasible and locally optimal. 

The proposed method has several strengths relative to other global optimization methods. It is agnostic of the forms of constraints in the problem; as long as we can query whether a sample $\tilde{\bx}$ is feasible to a constraint, we can embed the constraint into the \gls{mio} approximation. Once the constraints are learned using decision trees, the solution time of the resulting \gls{mio} approximation is low compared to solving the original global optimization problem. The proposed method can also be used to generate constraints from data which may not come from any known function, simulation or distribution. This allows us to simultaneously learn the physics of complex phenomena such as but not limited to social dynamical models or solutions of partial differential equations, and embed them into optimization problems.

In this work, we present our global optimization approach, implemented in our optimizer \gls{octhagon}, pronounced ``octagon''. We demonstrate its promise by considering global optimization problems with explicit nonlinear constraints. This allows us to quantify the performance of our method against existing global optimizers using available benchmarks. In addition, we approximate all nonlinear constraints in the benchmarks regardless of their efficient optimization-representability. The proposed method extends to mixed-integer-convex approaches where we embed efficiently-optimizable nonlinear constraints (e.g., quadratic, second order conic, log-sum-exponential constraints) into the \gls{mio} formulation directly, as long as these constraints are supported by the underlying solver. 

\subsection{Role of Machine Learning in Optimization}

The role of optimization in training \gls{ml} models is well known and studied. Recent review papers in the literature [\cite{Gambella2021,Sun2020}] survey the landscape of mathematical optimization and heuristic methods used for a variety of \gls{ml} applications. However, we are interested in the inverse of the above, and specifically how \gls{ml} can be used for the purpose of optimization, especially to solve problems that cannot naturally be posed as efficient optimization problems. 

There is precedent for using \gls{ml} methods to improve computational efficiency. A prominent example is the use of \gls{ml} to accelerate the simulation of nonlinear systems such as those in computational fluid dynamics~[\cite{Kochkov2021}], molecular dynamics~[\cite{Gastegger2017}] and quantum mechanics~[\cite{Morawietz2021}]. There has been some prior work using \gls{ml} to accelerate optimizations, e.g. using Bayesian optimization~[\cite{Frazier2018}] or neural networks~[\cite{Tagliarini1991}]. While these show that \gls{ml}-driven optimization is theoretically possible, the proposed methods are computationally expensive and not scalable for real-world problems. An interesting parallel use of \gls{ml} in optimization is in the interpretation of optimal solutions, where \gls{ml} is used to understand the optimal strategies (i.e. outcomes of all or subsets of decision variables) resulting from an optimization problem under different parameters~[\cite{Bertsimas2021}]. 

In this work, we use \gls{ml} to find optimal solutions to global optimization problems involving both explicit constraints with arbitrary mathematical primitives, and inexplicit black box functions. For this purpose, \gls{ml} is used for \emph{constraint learning} within two capacities. The first capacity is to accelerate optimizations over known models.  When models and/or constraints are known but their use is prohibitive, e.g. in the case of explicit but nonlinear and nonconvex constraints, learners are used to create surrogates that are more efficient for use in optimization. The second is in modeling. When data is available but models and/or constraints are black box, learners act as interpolants to the data, and to allow patterns in the data to be embedded in optimization.

Using \gls{ml} in optimization in both of these capacities requires that the approximating \gls{ml} models are optimization-representable. While many types of \gls{ml} models are efficiently queried and accurate, e.g., many types of neural networks and Gaussian processes, they cannot be embedded explicitly into structured optimization. Prior work has recognized the potential for using constraint learning approaches in optimization over data-driven constraints. Both \cite{Biggs2017} and \cite{Misic2020} use the prediction of tree ensembles as the objective function of optimization problems, given that a subset of tree features are decision variables. \cite{Maragno2021} go further and present a more general approach for data-driven optimization that leverages decision trees as well as other \gls{mio}-compatible \gls{ml} models such as support vector machines and neural networks. 

The aforementioned applications of decision trees in optimization are restricted in scope. \cite{Biggs2017} and \cite{Misic2020} limit their applications to optimization over data-driven objective functions, where decision trees are used to regress on a continuous quantity of interest. And while \cite{Maragno2021} use constraint learning for data-driven constraints, we use constraint learning to make approximations of intractable explicit and inexplicit constraints as well, where we have the capacity to sample the underlying constraints to generate data. Thus we propose a global optimization framework that can accommodate arbitrary explicit, inexplicit and data-driven constraints, leveraging decision trees in regression and classification settings. 

While it is possible to use other \gls{mio}-compatible \gls{ml} models for constraint learning in global optimization as proposed by \cite{Maragno2021}, we choose to rely on \gls{orth}s and \gls{octh}s since they are tunable, accurate and interpretable~[\cite{Bertsimas2019}]. In the following sections, we demonstrate that a global optimization method leveraging optimal decision trees makes significant progress in using \gls{ml} for both acceleration of optimizations and modeling, using the natural and intuitive \gls{mio} representation of trees. 

\subsection{Review of Decision Trees}

Decision trees is a popular predictive \gls{ml} method that partitions data hierarchically according to its features. A class label in a finite set of possible labels is assigned to each leaf node of the tree depending on the most common label of the data falling into the node. The optimization problem that is solved to produce a decision tree $T \in \mathbb{T}$ over known data $(\bx, \by)$ is the following:
\begin{equation*}
    \underset{T}{\mathrm{min}}~\mathrm{error}(T, \bx, \by) + c_p \cdot \mathrm{complexity}(T),
\end{equation*}
where $c_p$ is a complexity penalty parameter which attempts to strike a balance between the misclassification error over the test data and complexity (depth and breadth) of the tree. Once trained, decision trees are queried to predict the classes of test points with known features, but unknown class. 

Decision trees were pioneered by \cite{Breiman1984} with the advent of \gls{cart}. However, \gls{cart} is a top-down, greedy method of producing decision trees. Each split is only locally optimal since the splits are made recursively on the children of each new split starting from the root node. The ability of decision trees to explore the feature space has improved with the work of \cite{Bertsimas2017} on the \gls{oct}. \gls{oct}s leverage \gls{mio} and local search heuristics to reduce misclassification error relative to \gls{cart} without overfitting. Furthermore, \gls{oct}s are more \emph{interpretable}, since they can achieve similar misclassification error as trees generated by \gls{cart} with much less complexity.

\gls{octh}s generalize \gls{oct}s by allowing for hyperplane splits, i.e. splits in more than one feature at a time. An \gls{octh} can solve classification problems with higher accuracy and lower complexity than an \gls{oct} [\cite{Bertsimas2019}], and is more expressive in an optimization setting due to couplings of decision variables in nonlinear constraints. Thus, our method leverages \gls{octh}s exclusively to approximate constraints.

\gls{orth}s extend \gls{octh}s to regression problems, where the prediction of interest is continuous, i.e. $\tilde{y} \in \R$. Each leaf of an \gls{orth}, instead of containing a fixed class prediction, contains a continuous prediction $\tilde{y}$ as a linear regression over $\bx$ in the domain of the leaf. \gls{orth}s are particularly useful when approximating nonlinear objective functions. 

We rely on software from the company \gls{iai} in building, training and storing problem data in the form of \gls{octh}s and \gls{orth}s~[\cite{InterpretableAI}]. 

\section{Contributions}
\label{sec:gocontributions}

In this paper, we propose a global optimization approach that generalizes to explicit and inexplicit constraints and objective functions over bounded $\rm{dom}({\bx})$. Our specific contributions are as follows:
\begin{enumerate}
    \item We introduce an ensemble of methods for sampling constraints efficiently for the purpose of constraint learning. We leverage synergies of existing \gls{doe} techniques, but also devise a new \gls{knn} based sampling technique for sampling near-feasible points of explicit and inexplicit constraints. 
    \item We learn the feasible space of nonlinear objectives, inequalities and equalities using \gls{octh}s and \gls{orth}s. 
    \item We make \gls{mio} approximations of global optimization problems using the disjunctive representations of decision trees, and solve them using \gls{mio} solvers. 
    \item We devise a projected gradient descent method to check and repair the near-feasible, and near-optimal solutions from the \gls{mio} approximations. 
    \item We apply our method to a set of benchmark and real-world problems, and demonstrate its performance in finding global optima. 
\end{enumerate}

\subsection{Paper Structure.}

In Section~\ref{sec:method}, we detail our method, followed by a demonstrative example in Section~\ref{sec:demo}. In Section~\ref{sec:benchmarks}, we test our method on a number of benchmark problems from the literature, and compare our results with state-of-the-art global optimization tools such as BARON, IPOPT and CONOPT. In Section~\ref{sec:realworld}, we use our method to optimize two aerospace systems, one of which cannot be addressed via existing optimization tools. In Section~\ref{sec:godiscussion}, we discuss the results and avenues for future research. We conclude in Section~\ref{sec:goconclusion} by summarizing our findings and contributions. 

\section{Method}
\label{sec:method}

As aforementioned, our goal is to solve the global optimization problem approximately by making an \gls{octh} based \gls{mio} approximation, and then repairing the solution to be feasible and locally optimal. As an overview of this section, our method takes the following steps: 
\begin{enumerate}
    \item \textbf{Generate standard form problem:} In order to reduce the global optimization problem to a tractable \gls{mio} problem, we first restructure the global optimization problem in \eqref{eq:globalopt}. The linear constraints are passed directly to the \gls{mio} problem, while the nonlinear constraints are approximated in steps 2-6 below. If any variables involved in nonlinear constraints are unbounded from above and/or below, we attempt to compute bounds for the purpose of sampling. 
    \item \textbf{Sample and evaluate nonlinear constraints:} The data used in training is important for the accuracy of \gls{ml} models. For accurate \gls{octh} approximations of nonlinear constraints, we use fast heuristics and \gls{doe} methods to sample variables over $\rm{dom}(\bx)$. We evaluate each constraint over the samples, and resample to find additional points near the constraint boundary for local approximation refinement. 
    \item \textbf{Train decision trees over constraint data}: The feasibility space of each constraint is classified and approximated by an \gls{octh}. If the objective function is nonlinear, it is regressed and approximated via an \gls{orth}.
    \item \textbf{Generate \gls{mi} approximation:} The decision paths and hyperplane splits are extracted from the trees, and used to formulate efficient \gls{mio} approximations of the nonlinear constraints using disjunctions.  
    \item \textbf{Solve \gls{mio} approximation:} The resulting \gls{mio} problem is optimized using commercial solvers to get an approximate solution.  
    \item \textbf{Check and repair solution:} The \gls{mio} problem approximates the global optimization problem, so the optimum is likely to be near-optimal and near-feasible. We evaluate the feasibility of each nonlinear constraint, and compute the gradients of the objective and nonlinear constraints using automatic differentiation. In case of suboptimality or infeasibility, we perform a number of projected gradient descent steps to repair the solution, so that it is feasible and locally optimal. 
\end{enumerate}

We describe the steps in greater detail in Sections~\ref{sec:stdform} through \ref{sec:repair}. A step-by-step demonstration of the method, as implemented in our optimizer \gls{octhagon}, can be found in Section~\ref{sec:demo}.

\subsection{Standard Form Problem.}
\label{sec:stdform}

We restructure the global optimization problem posed in \eqref{eq:globalopt} by separating the linear and nonlinear constraints. The linear constraints are passed directly into a \gls{mio} model, while the nonlinear constraints are stored for approximation. If constraints are black box, they are assumed to be nonlinear as well. This restructured problem is shown in \eqref{eq:stdform}, and referred to as the standard form. Note that the standard form allows for both nonlinear inequalities and equalities.
\begin{align}
    \begin{split}
        \underset{x}{\text{min}} &~~ f(\bx) \\
        \text{s.t.} &~~ g_i(\bx) \geq 0,~ i \in I, \\ 
        &~~ h_j(\bx) = 0,~j \in J, \\ 
        &~~\mathbf{Ax} \geq \mathbf{b},~\mathbf{Cx} = \mathbf{d}, \\
        &~~x_k \in [\underline{x}_k, \overline{x}_k],~k \in [n]. 
        \label{eq:stdform}
    \end{split}
\end{align}
\vspace{0.25in}
\subsubsection{Variable Outer-Bounding.}

The proposed method requires boundedness of decision variables $\bx$ in each approximated constraint so that we can sample $\rm{dom}(\bx)$ for constraint evaluation. When bounds are missing for any variable $x_k$ in a nonlinear constraint, we pose the following optimization problem over the linear constraints only. 
\begin{align}
    \begin{split}
         \underset{x}{\text{min/max}} &~~ x_k \\
        \text{s.t.} &~~\mathbf{Ax} \geq \mathbf{b},~\mathbf{Cx} = \mathbf{d} \\
        &~~x_i \in [\underline{x}_i, \overline{x}_i],~i \subseteq [n].
    \end{split}
\end{align} 

The solution to this problem is the absolute largest range $[\underline{x}_k, \overline{x}_k]$ that satisfies all linear constraints as well as bounds on $x_i$, for those indices $i$ for which $x_i$ is bounded. We can also solve the above optimization problem to tighten bounds on variables with existing bounds. Tighter bounds can significantly improve solution quality and time by improving the quality of \gls{ml} approximations.

\subsection{Sampling and Evaluation of Nonlinear Constraints.}
\label{sec:sampling}

For the purpose of constraint learning, we require data over variables and corresponding left-hand-side values of nonlinear constraints. The importance of the quality of data for the accuracy of machine learning tasks is well known and studied since the 1990's~[\cite{Cortes1995}]. Thus, the distribution of data points used for constraint learning is critical. The samples over $\rm{dom}(\bx)$ should be sufficiently space-filling so that the behavior of each constraint is captured over the whole $\rm{dom}(\bx)$. In addition, we require sufficient concentration of points near the constraint boundary so that learners are adequately trained to predict the feasibility of near-feasible points. 

To achieve both of these objectives, we take a disciplined approach to sampling, and generate data over $\rm{dom}(\bx)$ for each constraint in several stages. Note that the sampling and evaluation steps in the following subsections are performed constraintwise. 

\subsubsection{Boundary Sampling.}

We first sample the corners of the $\bx$ hypercube for the constraint, defined by $x_k \in [\underline{x}_k, \overline{x}_k],~k \subseteq [n]$, in an effort to capture extremal points. We call this boundary sampling. This is combinatorial in the number of variables in each nonlinear constraint; a constraint with $p$ bounded variables would require $2^p$ samples. In practice, we sample a limited combination of corner points, depending on the number of variables in the constraint. 

\subsubsection{Optimal Latin Hypercube Sampling.}

Next, we implement \gls{olh} sampling over the $\bx$ hypercube. There is a wealth of literature starting with \cite{McKay1979} that demonstrates the strength of \gls{lh} sampling versus other methods for \gls{doe}. However, \gls{lh} sampling is not in general a \emph{maximum entropy sampling scheme}~[\cite{Shewry1987}], i.e. the samples from \gls{lh}s do not optimize information gained about the underlying system. \gls{olh} sampling is the entropy maximizing variant of \gls{lh} sampling for a uniform prior, where our entropy function is the pairwise Euclidian distances between sample points~[\cite{Bates2003}]. The uniform prior assumption is logical since we do not have or require an initial guess for where in the $\bx$ hypercube the optimal solution will land, and the constraints are treated as black boxes. 

\gls{olh} sampling, unlike standard \gls{lh}s, is space-filling and thus useful for learning the global behavior of constraints using \gls{ml} models. In practice, \gls{olh} generation is time-consuming and impractical. Instead, we use an efficient heuristic proposed by \cite{Bates2004}, which uses a permutation genetic algorithm to find near-optimal solutions to the \gls{olh} problem with low computational cost. We terminate the genetic algorithm prematurely in our optimization scheme, since samples are not required to be optimally distributed. 

\subsubsection{Constraint Evaluation.}

We use the samples to either compute the left-hand-side of the constraint, or the feasibility of the constraint if the left-hand-side is not available. If the constraint is an equality $h_j(\bx) = 0$, we relax it and treat it as an inequality $h_j(\bx) \geq 0$ until Section~\ref{sec:mi_constraints}. The result is a $\{0, 1\}^n$ feasibility vector corresponding to each of the $n$ samples, defining the classes for the classification problem. 

If desired, assuming that constraints use a common set of samples, it is possible to lump the feasibility of a set of inequality constraints by taking the row-wise minimum of their joint feasibility over the same data. This can reduce the model complexity, but we currently do not consider this in our method. 

\subsubsection{$k$NN Quasi-Newton Sampling.}

The previous sampling methods achieve a space-filling distribution of samples in $\rm{dom}(\bx)$ to enable approximating \gls{octh}s to learn the feasibility of each constraint in a global sense. We still require sufficient concentration of points near the constraint boundary, i.e. points $\tilde{\bx}_i$ so that $g(\tilde{\bx_i}) \approx 0$, so our \gls{octh} models are trained to classify such near-feasible points accurately. 

Assuming that the first stage sampling and evaluation has found at least one feasible point to the constraint, in this step, we attempt to sample near the constraint boundary using a method we've developed called \gls{knn} quasi-Newton sampling. The method hinges on using \gls{knn} to generate near-feasible neighborhoods for the constraint over previous data $(\tilde{\bx}, \tilde{\by})$, and using approximate gradients in these neighborhoods to find new near-feasible samples $\tilde{\bu}$, with vanishing $g(\tilde{u}_i) = \epsilon \rightarrow 0$. We present the method in Algorithm~\ref{alg:knnsampling}.\\

\begin{algorithm}[H]
\caption{\gls{knn} quasi-Newton sampling}
\label{alg:knnsampling}
\KwResult{Sample points near the feasibility boundary of constraints.}
    Find $k = p+1$ nearest neighbors: $\mathbf{\xi} = [k\mathrm{NN}(\tilde{\bx}_{i},~ \tilde{\bx}),~\forall i \in [n]$]\;
    Classify feasibility \gls{knn} patches: $\mathbf{\phi} \in \{\rm{feasible, infeasible, mixed}\}^n$\;
    Initialize new sample container $\tilde{\bu} = []$. \;
    \For{$i \in [n]$:}{
        \If{$\phi_i = \rm{mixed}$ and $\tilde{\bx}_i$ \rm{infeasible}} {
            \For{$j \in [k]$}{
            \If{$\tilde{\bx}_{\xi_{i,j}}$ \rm{feasible}} {
                Augment $\tilde{\bu}$: secant method($\tilde{\bx}_{i},~\tilde{y}_{i},~\tilde{\bx}_{\xi_{i,j}},~\tilde{y}_{\xi_{i,j}}$)
            }
        }
        }
    }
\end{algorithm}

The method is described as follows. Starting from space-filling data $(\tilde{\bx},\tilde{\by})$ where $\tilde{y}_{i} = g(\tilde{\bx}_{i})$, we find the $k$-nearest points for each sampled point $\tilde{\bx}_i$ in the 0-1 normalized $\bx$ hypercube. In our particular implementation, we use $k=p+1$, where $p$ is the number of variables in constraint $g(\bx) \geq 0$. For each \gls{knn} cluster with index $i$ centered at $\tilde{\bx}_i$ with $k-1$ neighbor indices $\xi_i$, we determine if all sample points are feasible, all points are infeasible, or points are mixed-feasibility. 

In each cluster with mixed-feasibility points, we perform the secant method between points of opposing feasibility. The secant method is an approximate root finding algorithm defined by the following recurrence relation
\begin{equation}
    \tilde{\bx}_k = \tilde{\bx}_j - \tilde{y}_j\frac{\tilde{\bx}_j - \tilde{\bx}_{i}}{\tilde{y}_j - \tilde{y}_i},
\end{equation}
where $\tilde{\bx}_i$ and $\tilde{\bx}_j$ are points of opposing feasibility in the same mixed-feasibility neighborhood, and $\tilde{\bx}_k$ is a new candidate root. The secant method thus allows us to efficiently generate roots $\tilde{\bx}_k$ that would be expected to be near the constraint boundary, using combinations of points $\tilde{\bx}_i$ and $\tilde{\bx}_j$ from the space-filling \gls{olh} samples. We ensure that each pair of \gls{knn}-adjacent points on the constraint boundary results in only one new point, by only sampling within mixed-feasibility \gls{knn} cells if their centroid is infeasible, and then only sampling between the infeasible centroid and surrounding feasible points in the \gls{knn} cell. 

Once we have performed the \gls{knn} sampling process and have new samples $\tilde{\bu}$, we evaluate the left-hand-side $g(\bx)$ over the samples and add them to data $(\tilde{\bx},\tilde{\by})$ before proceeding to the tree training step. 

\subsection{Decision Tree Training.}

We use trees to approximate the nonlinear constraints in our global optimization problem due to their \gls{mio} representability, which we will demonstrate in Section~\ref{sec:mi_constraints}. We use software from the company \gls{iai} in building, training and storing problem data in the form of \gls{octh}s and \gls{orth}s [\cite{InterpretableAI}]. We train trees exclusively with hyperplane splits due to their higher approximation accuracy and lower tree complexity.

The trees are trained on all available data instead of a subset of the data as would be expected in traditional \gls{ml}. In addition, we penalize tree complexity very little. This is because our data is noise-free, and approximation accuracy is important in the global optimization setting. In the case where the constraints are generated on noisy data, we would allow for the splitting of data into training and test sets, and cross-validate over a range of parameters.

We use the base \gls{octh} and \gls{orth} parameters in Table~\ref{tab:treeparams} within \gls{iai} when initializing constraint learning instances. These parameters are used for all computational benchmarks throughout the paper unless stated otherwise. The parameters have been chosen to balance tree accuracy with tree complexity and associated computational cost, and may be tuned by users as they find necessary. 
\begin{center}
    \begin{table}[h!]
        \centering
            \begin{tabular}{|c|c|c|}
                \hline
                Parameter & \gls{octh} & \gls{orth} \\
                \hline
                Hyperplane sparsity & All & All \\
                Regression sparsity & - & All \\
                Max depth & 5 & 5 \\
                Complexity factor & $10^{-6}$ & $10^{-6}$\\
                Minbucket & 0.01 & 0.02 \\
                Random tree restarts & 10 & 10 \\
                Hyperplane restarts & 5 & 5 \\
                \hline
            \end{tabular}
            \caption{Parameters for base decision trees in constraint learning.}
            \label{tab:treeparams}
\end{table}
\end{center}

Our training loss function for \gls{octh}s is misclassification error. 
If a tree is a function that maps feature inputs into classes ($T: \bx \xrightarrow{}  y$), the misclassification error is simply the weighted proportion of samples that are misclassified by the tree, where $\I$ is the indicator function and $w_i$ are the sample weights. An exact classifier would have a misclassification error of 0. 
\begin{equation*}
    \text{misclassification error} = \frac{1}{n}\frac{\sum_{i=1}^{n} w_i \cdot \I(T(\bx_i) \neq y_i)}{\sum_{i=1}^{n} w_i}.
\end{equation*}

For \gls{orth}s used to approximate objective functions, we use $1-\rm{R}^2$ as the loss function, where $\rm{R}^2$ is the coefficient of determination. An exact regressor would have a $1-\rm{R}^2$ value of 0. 
\begin{equation*}
    1-\mathrm{R}^2 = \frac{\sum_{i=1}^n(T(\bx_i) - y_i)^2}{\sum_{i=1}^n(T(\bx_i) - \bar{y})^2},~\mathrm{where}~\bar{y} = \frac{1}{n} \sum_{i=1}^n y_i.
\end{equation*}

\subsection{MI Approximation.}
\label{sec:mi_constraints}

From this section forward, we recognize that the global optimization problem is approximated constraint-wise, and introduce indices $i \in I$ and $j \in J$ for the inequality and equality constraints respectively. Having classified the feasible space of nonlinear inequalities $g_i(\bx) \geq 0,~i \in I$ and relaxed nonlinear equalities $h_j(\bx) \geq 0,~j \in J$ using \gls{octh}s, we retighten equalities to $h_j(\bx) = 0,~j \in J$, and pose the feasible $\bx$-domains of each tree as unions of polyhedra. In this section, we define mathematically the set of disjunctive \gls{mi}-linear constraints that represent the trees exactly. 

\subsubsection{Nonlinear Inequalities.}

The tree $T_{i}$ that classifies the feasible set of nonlinear inequality $g_i(\bx) \geq 0$ has a set of leaves $L_{i}$, where a subset of leaves $L_{i,1} \subset L_i$ are classified feasible (where the indicator function $\I(g_i(\bx) \geq 0) = 1$) and $L_{i,0} \subset L_i$ are classified infeasible ($\I(g_i(\bx) \geq 0) = 0$). The decision path to each leaf defines a set of separating hyperplanes, $H_{i,l}$, where $H_{i,l,-}$ and $H_{i,l,+}$ are the set of leftward (less-than) and rightward (greater-than) splits required to reach leaf $l$ respectively. The feasible polyhedron of tree $T_i$ at feasible leaf $l \in L_{i,1}$ is thus defined as 
\begin{equation}
    \label{eq:leaf}
    \bbp_{i,l} = \{\bx: \balpha_h^{\top} \bx \leq \beta_h, ~\forall~h \in H_{i,l,-}~;~ \balpha_h^{\top} \bx \geq \beta_h, ~\forall~ h \in H_{i,l,+}\}.
\end{equation}

The feasible set of $\bx$ over constraint $g_i(\bx) \geq 0$ is approximated by the union of the feasible polyhedra in \eqref{eq:leaf}. More formally, 
\begin{equation}
    \bx \in \bigcup_{l \in L_{i,1}} \bbp_{i,l}.
\end{equation}
This union-of-polyhedra representation can described by a set of disjunctive constraints involving a big-M formulation. \cite{Vielma2015} describes many such ``projected" formulations; the specific disjunctive representation of \gls{octh}s approximating nonlinear inequalities is as follows:
\begin{equation}
\bx \in \bigcup_{l \in L_{i,1}} \bbp_{i,l} \iff
\begin{cases}
    \begin{aligned}
    \label{eq:disjunctive_M}
    &\{~\balpha_h^{\top} \bx \leq \beta_h + M(1-z_{i,l}), ~\forall~h \in H_{i,l,-}~; \\
    &~\beta_h \leq \balpha_h^{\top} \bx  + M(1-z_{i,l}) , ~\forall~ h \in H_{i,l,+}\},~ \forall~l \in L_{i,1}, \\
    ~&\sum_{l \in L_{i,1}} z_{i,l} = 1, \\
    &~ z_{i,l} \in \{0, 1\},~l \in L_{i,1},  \\
    &~M > |\beta_h|,~M > \underset{\rm{dom}(\bx)}{\max}~|\balpha_h^{\top} \bx|,~\forall~h \in H_{i,l},~l \in L_{i,1}.
    \end{aligned}
\end{cases}
\end{equation}

Membership of $\bx$ in polyhedron $\bbp_{i,l}$ is defined by binary variable $z_{i,l}$. The constraint $\sum_{l \in L_{i,1}} z_{i,l} = 1$ ensures that $\bx$ is in exactly one feasible polyhedron. However, the formulation above requires knowing the value of $M$ with sufficient accuracy, which can be difficult in practice. The value of $M$ is important; too small an $M$ means that the constraint is insufficiently enforced, and too large an $M$ can cause numerical issues. Knowing M to a sufficient tolerance can require solving the inner maximization in \eqref{eq:disjunctive_M} over $\rm{dom}(\bx)$, and even declaring a separate $M_h$ for each separating hyperplane $h \in H_{i,l}$. 

Alternatively, we derive a representation that completely avoids the need to compute big-M values, since we restrict ourselves to $\bx \in [\underline{\bx}, \overline{\bx}]$. The tradeoff is that we require the addition of auxiliary variables $\by_l \in \R^{p_i}$, for each leaf $~l \in L_{i,1}$, where $p_i$ is the dimension of variables in constraint $i$. We present the big-M free representation of \gls{octh}s used to approximate nonlinear inequalities in \eqref{eq:disjunctive}. The formulation is an application of basic extended disjunctive formulations for defining unions of polyhedra, as detailed by \cite{Vielma2015}.
\begin{align}
    \bx \in \bigcup_{l \in L_{i,1}} \bbp_{i,l} \iff
    \begin{cases}
        \begin{split}
        \label{eq:disjunctive}
        &\{~\balpha_h^{\top} \by_l \leq \beta_h z_{i,l}, ~\forall~h \in H_{i,l,-}~; \\
        &~\beta_h z_{i,l} \leq \balpha_h^{\top} \by_l, ~\forall~ h \in H_{i,l,+}\}~\forall~l \in L_{i,1}, \\
        &~\by_l \in [\underline{\bx} z_{i,l}, \overline{\bx} z_{i,l}],~l \in L_{i,1}, \\
        &\sum_{l \in L_{i,1}} \by_l = \bx, \\
        ~&\sum_{l \in L_{i,1}} z_{i,l} = 1, \\
        &~ z_{i,l} \in \{0, 1\},~l \in L_{i,1}.  \\
        \end{split}
    \end{cases}
\end{align}

Just as the big-M formulation, whether or not $\bx$ lies in polyhedron $\bbp_{i,l}$ is defined by binary variable $z_{i,l} \in \{0, 1\}$. If $\bx$ is in $\bbp_{i,l}$, then $\bx = \by_l$. If not, $\by_l = \mathbf{0}$. Thus $\bx$ can only lie in the leaves of $T_i$ that are classified feasible. 

Notably, formulation~\eqref{eq:disjunctive} is \emph{locally ideal}, i.e. its continuous relaxation has at least one basic feasible solution, and all its basic feasible solutions are integral in $\bz_i$~[\cite{Vielma2015}]. This confers computational advantages in optimization over such disjunctions compared to its big-M variant. Since disjunctive formulation \eqref{eq:disjunctive} is tractable and big-M free, we implement it in \gls{octhagon}. 

\subsubsection{Nonlinear Equalities.}

Nonlinear equalities can also be approximated by \gls{octh}s. To do so, we simply relax $h_j(\bx) = 0$ to $h_j(\bx) \geq 0$ and fit an \gls{octh} $T_j$ to the feasible set of this constraint, with polyhedra $\bbp_{j,l}$, where $l$ can lie in feasible leaves $L_{j,1}$ and infeasible leaves $L_{j,0}$. The feasible set of the original equality must be represented by the union of the polyhedral faces between the feasible and infeasible leaves. It is critical to note however that this is not equivalent to the union of polyhedral faces, $\bx \in \bigcup_{l \in L_{j}} \text{faces} (\bbp_{j,l})$, since some of the faces separate two feasible spaces from each other, and thus would not be valid constraint boundaries. We are only interested in polyhedral faces that separate feasible polyhedra from infeasible polyhedra, where $h_j(\bx) \geq 0$ and $h_j(\bx) \leq 0$. 
Therefore the approximate equality is the union of intersections of all permutations of a feasible polyhedron with an infeasible polyhedron,
\begin{equation}
    \bx \in \bigcup_{l_0 \in L_{j,0},~l_1 \in L_{j,1}} \{\bbp_{j,l_0} \cap \bbp_{j,l_1}\}. 
\end{equation}

To ensure that $\bx$ lies on a face between a feasible and an infeasible polyhedron, we allocate a binary variable $z_{j,l}$ for each leaf $l \in L_j$. We make sure that $\bx$ lies in exactly one feasible and one infeasible polyhedron by having exactly two non-zero $z_{j,l}$'s, one in a feasible leaf $l \in L_{j,1}$ and the other in an infeasible leaf $l \in L_{j,0}$. Thus we  represent the approximate equality as the following set of disjunctive big-M constraints, where $L_j = \{L_{j,1} \cup L_{j,0}\}$ are the combined set of feasible and infeasible leaves of tree $T_j$.
\begin{align}
\bx \in \bigcup_{\substack{l_0 \in L_{j,0},\\l_1 \in L_{j,1}}} \{\bbp_{j,l_0} \cap \bbp_{j,l_1}\} \iff
\begin{cases}
    \begin{split}
    \label{eq:disjunctive_eq_M}
    &\{~\balpha_h^{\top} \bx \leq \beta_h + M(1-z_{j,l}), ~\forall~h \in H_{j,l,-}~; \\
    &~\beta_h \leq \balpha_h^{\top} \bx  + M(1-z_{j,l}) , ~\forall~ h \in H_{j,l,+}\},~\forall l \in L_j, \\
    ~&\sum_{l \in L_{j,0}} z_{j,l} = 1,~~\sum_{l \in L_{j,1}} z_{j,l} = 1, \\
    &~ z_{j,l} \in \{0, 1\},~~~l \in L_j.  
    \end{split}
\end{cases}
\end{align}
This guarantees that $\bx$ falls on a polyhedral face that separates a feasible and infeasible polyhedron, thus approximating $h_j(\bx) = 0$. As we have done for nonlinear inequalities, we can come up with an equivalent big-M-free formulation as follows, and implement it in \gls{octhagon}.
\begin{align}
    \bx \in \bigcup_{\substack{l_0 \in L_{j,0},\\l_1 \in L_{j,1}}} \{\bbp_{j,l_0} \cap \bbp_{j,l_1}\} \iff
    \begin{cases}
        \begin{split}
        \label{eq:disjunctive_eq}
        &\{~\balpha_h^{\top} \by_l \leq \beta_h z_{i,l}, ~\forall~h \in H_{i,l,-}~; \\
        &~\beta_h z_{i,l} \leq \balpha_h^{\top} \by_l, ~\forall~ h \in H_{i,l,+}\},~ \forall~l \in L_j, \\
        &~\by_l \in [\underline{\bx} z_{i,l}, \overline{\bx} z_{i,l}],~l \in L_j, \\
        &\sum_{l \in L_{i,1}} \by_l = \bx,~~\sum_{l \in L_{i,0}} \by_l = \bx, \\
        ~&\sum_{l \in L_{i,1}} z_{i,l} = 1,~~\sum_{l \in L_{i,0}} z_{i,l} = 1, \\
        &~ z_{i,l} \in \{0, 1\},~~~l \in L_j.  \\
        \end{split}
    \end{cases}
\end{align}

Note that nonlinear equalities pose the greatest challenge for any global optimization method, since the $\epsilon$-feasible space of equalities is restrictive. 

\subsubsection{Nonlinear Objectives.}
\label{sec:nlobjectives}

We treat nonlinear objectives $f(\bx)$ differently than constraints. Constraints are represented well by classifiers because constraints partition the space of $\bx$ into feasible and infeasible classes. Nonlinear objectives however are continuous with respect to $\bx$, and are thus better approximated by regressors. To approximate a nonlinear objective function $f(\bx)$, we train an \gls{orth} on sample data $\{\tilde{\bx}_i, f(\tilde{\bx}_i)\}_{i=1}^n$, and replace the nonlinear objective with the auxiliary variable $f^*$. We lower bound the value of $f^*$ using the disjunctive constraints derived from the \gls{orth}, thus approximating the original objective function. 

We can apply the same logic to constraints of the form  $\mathbf{a}^{\top} \bx + b \geq g(\bx)$, where the left-hand-side is affine  and separable from the nonlinear component $g(\bx)$. Since $\mathbf{a}^{\top} \bx + b$ is linear and \gls{mio}-compatible, we instead train an \gls{orth} on sample data $\{\tilde{\bx}_i, g(\tilde{\bx}_i)\}_{i=1}^n$, and make sure that $\mathbf{a}^{\top} \bx + b$ is lower bounded by the approximating \gls{orth}. It is the choice of the user whether or not to use \gls{octh}s or \gls{orth}s to approximate separable constraints, but in general an \gls{orth} is more accurate in these cases. All problems addressed in Section~\ref{sec:benchmarks} treat constraints as non-separable, and use classifiers to approximate them. To solve the satellite scheduling problem in Section~\ref{sec:oos}, we take advantage of this separability and choose to train \gls{orth}s instead.

Since an \gls{orth} is an \gls{octh} with additional regressors added to each leaf, the disjunctive constraints in \eqref{eq:disjunctive} and \eqref{eq:disjunctive_eq} apply with minor modifications described as follows. $L_f$ is the set of leaves of the approximating \gls{orth}; assuming that $f(\bx)$ can be evaluated on $\rm{dom}(\bx)$, all leaves $l \in L_f$ of the \gls{orth} can feasibly contain $\bx$, meaning that the disjunctions are applied to all leaves instead of a subset of the leaves of the tree. Each leaf $l \in L_f$ has a set of separating hyperplanes that is described by its decision path, as well as an additional separating hyperplane described by the regressor in each leaf.

For objectives and separable inequalities, instead of using the regressor within each leaf of the \gls{orth} directly, we run a secondary linear regression problem on the points within each leaf to find the tightest lower bounding hyperplane on the data. This allows us to have an approximate relaxation of the constraint or objective function, and tighten the relaxation later via solution repair in Section~\ref{sec:repair}. 

\subsection{Solution of MIO Approximation.}
\label{sec:mio}

Having represented the feasible space of inequality and equality constraints as a unions of polyhedra, we have the following final problem. 
\begin{align}
    \begin{split}
        \underset{x}{\text{min}} &~~ f^*\\
        \text{s.t.} &~~ f^*, \bx \in \bigcup_{l \in L_f} \bbp_{i,l}, \\
        &~~\bx \in \bigcup_{l \in L_{i,1}} \bbp_{i,l},~ \forall~i \in I, \\ 
        &~~ \bx \in \bigcup_{l_0 \in L_{j,0},~l_1 \in L_{j,1}} \{\bbp_{j,l_0} \cap \bbp_{j,l_1}\},~\forall j \in J, \\
        &~~\mathbf{Ax} \geq \mathbf{b},~\mathbf{Cx} = \mathbf{d}, \\
        &~~ x_k \in [\underline{x}_k, \overline{x}_k],~k \in [n].
        \label{eq:finalproblem}
    \end{split}
\end{align}

This is a \gls{milo} that can be efficiently solved using branch-and-bound methods. We use CPLEX for this purpose, since it is available free of charge to solve small scale \gls{milo} instances. 

\subsection{Solution Checking and Repair.}
\label{sec:repair}

The optimum obtained in Section~\ref{sec:mio} is likely to be near-optimal and near-feasible to the original global optimization problem, since the \gls{mio} is approximate. To repair the solution in case of suboptimality or infeasibility, we devise and present a local search procedure based on \gls{pgd}. \gls{pgd} is a method for constrained gradient descent that is reliable, scalable and fast for the local optimization required to restore feasibility and optimality to approximate solutions. It relies on using gradients of the constraints and objective to simultaneously reduce constraint violation (by projecting $\bx^*$ onto the feasible space of $\bx$) and the objective function value. Our particular implementation of \gls{pgd} solves a series of gradient-driven \gls{mio} problems to do so.

To obtain the gradients of explicit and inexplicit constraints, we leverage \gls{ad}, and specifically \emph{forward mode \gls{ad}}. Forward mode \gls{ad} looks at the fundamental mathematical operations involved in evaluating the constraint functions, and thus computes the gradient of each constraint exactly at any solution $\bx^*$~[\cite{Verma2000}]. Unlike finite differentiation, \gls{ad} does not require additional function evaluations or discretization, and unlike symbolic differentiation, it doesn't require the constraints to be explicit. 

The proposed \gls{pgd} method begins by first evaluating the objective and all constraints at $\bx^*$, the last known optimum, as well as their gradients. The disjunctive approximations of nonlinear inequality constraints are replaced by linear approximators based on the local constraint gradient, depending on the feasibility of each constraint:  
\begin{align}
    g_i(\bx) \geq 0 \rightarrow
    \begin{cases}
        \begin{split}
            \nabla g_i(\bx^*)^{\top} \bd + g_i(\bx^*) \geq 0, &~\mathrm{if}~ g_i(\bx^*) \geq 0, \\[1ex]
            \nabla g_i(\bx^*)^{\top} \bd + g_i(\bx^*) + \lambda_i \geq 0, & ~\mathrm{if}~ g_i(\bx^*) \leq 0, \\
        \end{split}
        \label{eq:pgdinequality}
    \end{cases}
\end{align}
where $\bd \in \R^n$ is the descent direction, and $\lambda_i \in \R^+$ is an inequality relaxation variable. Similarly, we replace the \gls{mi} approximations of equalities with their local linear approximators, but always include relaxation variables regardless of the level of infeasibility of the constraints, as shown in \eqref{eq:pgdequality}.
\begin{align}
    h_j(\bx) = 0 \rightarrow
    \begin{cases}
        \begin{split}
            \nabla h_j(\bx^*)^{\top} \bd + h_j(\bx^*) + \mu_j &\geq 0, \\[1ex]
            \nabla h_j(\bx^*)^{\top} \bd + h_j(\bx^*) &\leq \mu_j, \\
        \end{split}
        \label{eq:pgdequality}
    \end{cases}
\end{align}
where $\mu_j \in \R^+$ is an equality relaxation variable. This relaxation is for two reasons. The first is that, in presence of equalities, the local \gls{pgd} step may be infeasible due to conflicting equality constraints. The second is that each \gls{pgd} step will involve solving a quadratic program, which can only be solved to given numerical precision. This precision, while low, is non-zero. 

Thus we introduce a constraint tightness tolerance parameter $\phi$, and say that an inequality $g_i(\bx) \geq 0$ is feasible at $\bx^*$ if $g_i(\bx^*) \geq -\phi$. If all inequality constraints are feasible to tolerance, relaxation variables $\mathbb{\lambda}$ are only required on the inequalities where $0 \geq g_i(\bx^*) \geq -\phi,~i \in I$, by the condition in \eqref{eq:pgdinequality}. In that case, we perform a simple gradient descent step. This involves solving the quadratic optimization problem in \eqref{eq:simpledescent}, where $\gamma$ is the infeasibility penalty coefficient, $\alpha$ is the step size within a 0-1 normalized $\bx$ hypercube, $r$ is the step size decay rate, $t$ is the current \gls{pgd} iteration and $T$ is the maximum number of iterations.  
\begin{equation}
    \begin{aligned}
            \underset{\bx, \bd, \lambda, \mu}{\text{min}} &~~ \nabla f(\bx^*)^{\top} \bd + \gamma (||\mathbb{\lambda}||_2^2 + ||\mathbb{\mu}||_2^2) \\
            \text{s.t.} &~~ \bx = \bx^* + \bd, \\
            &~~ \Bigg|\Bigg| \frac{\bd}{\overline{\bx} - \underline{\bx}} \Bigg|\Bigg|^2_2 \leq \alpha \rm{exp} \Big( \frac{-rt}{T} \Big), \\
            &~~ \left\{\begin{array}{lr} 
                \nabla g_i(\bx^*)^{\top} \bd + g_i(\bx^*) \geq 0, &~\mathrm{if}~ g_i(\bx^*) \geq 0 \\
                \nabla g_i(\bx^*)^{\top} \bd + g_i(\bx^*) + \lambda_i \geq 0, & ~\mathrm{if}~ -\phi \leq g_i(\bx^*) \leq 0
                \end{array}\right\},~\forall i \in I, \\
            &~~ \left\{\begin{array}{lr}
                \nabla h_j(\bx^*)^{\top} \bd + h_j(\bx^*) + \mu_j \geq 0, \\
                \nabla h_j(\bx^*)^{\top} \bd + h_j(\bx^*) \leq \mu_j,
                \end{array}\right\},~\forall j \in J, \\
            &~~\mathbf{Ax} \geq \mathbf{b},~\mathbf{Cx} = \mathbf{d}, \\
            &~~ x_k \in [\underline{x}_k, \overline{x}_k],~\forall k \in [n] \\
            &~~ \left\{\begin{array}{lr} 
                \lambda_i = 0, &~\mathrm{if}~ g_i(\bx^*) \geq 0 \\
                \lambda_i \geq 0, & ~\mathrm{if}~ g_i(\bx^*) \leq 0
                \end{array}\right\},~\forall i \in I, \\
            &~~ \mu_i \in \R_+,~j \in J.
        \label{eq:simpledescent}
    \end{aligned}
\end{equation}
We exponentially decrease the allowed step size $\bd$ as defined in \eqref{eq:simpledescent}, to aid convergence and break cycles that may result. 

If the current solution $\bx^*$ is infeasible beyond tolerance to any constraints, we take a projection-and-descent step. This modifies the objective and first two constraints in ~\eqref{eq:simpledescent} by removing the step size constraint on $\bd$, and augmenting the objective function with a projection distance penalty with $\beta$ as a parameter, as shown in \eqref{eq:projecteddescent}:
\begin{align}
    \begin{split}
        \underset{\bx, \bd, \lambda, \mu}{\text{min}} &~~ \nabla f(\bx^*)^{\top} \bd + \beta \Bigg|\Bigg| \frac{\bd}{\overline{\bx} - \underline{\bx}} \Bigg|\Bigg|^2_2 + \gamma (||\mathbb{\lambda}||_2^2 + ||\mathbb{\mu}||_2^2) \\
        \text{s.t.} &~~ \bx = \bx^* + \bd, \\
        \vdots
        \label{eq:projecteddescent}
    \end{split}
\end{align}
This quadratic optimization problem approximates the closest feasible projection of $\bx$ onto the feasible space of nonlinear constraints. 

The gradient and projected gradient steps defined above require knowing the maximum range on all variables, $\overline{\bx} - \underline{\bx}$. If this range is not provided for variable $x_k$, then we assume $\overline{x}_k - \underline{x}_k = \max(\overline{\bx}) - \min(\underline{\bx})$. The convergence of the \gls{pgd} is much stronger with user-provided bounds however. We repeat the above \gls{pgd} steps on new incumbent solutions until the final two solutions are feasible to all constraints, and the improvement in original objective function $f(\bx)$ is less than absolute tolerance $\epsilon$. 

The \gls{pgd} algorithm introduces many parameters, whose default values are defined in Table~\ref{tab:pgdparams}. While this adds additional complexity to the solution procedure, the descent procedure is intuitive to tune, and the current implementation warns the user in case parameters require examination. In addition, the parameters are applied to 0-1 normalized quantities over the $\bx$ hypercube wherever possible. For all examples in this paper, the default \gls{pgd} parameters from Table~\ref{tab:pgdparams} apply unless stated otherwise. 

\begin{center}
    \begin{table}[h!]
        \centering
            \begin{tabular}{|c|c|c|}
                \hline
                Parameter & Description & Value \\
                \hline
                $\gamma$ & Infeasibility penalty & $10^{6}$ \\
                $\beta$ & Step penalty & $10^4$ \\
                $\alpha$ & Step size & $10^{-3}$ \\
                $r$ & Decay rate & 2 \\
                $T$ & Maximum iterations & 100 \\
                $\epsilon$ & Absolute tolerance & $10^{-4}$ \\
                $\phi$ & Tightness tolerance & $10^{-8}$ \\
                \hline
            \end{tabular}
            \caption{Parameters for PGD repair procedure.}
            \label{tab:pgdparams}
\end{table}
\end{center}

\section{Demonstrative Example}
\label{sec:demo}

Consider the following modified mixed-integer nonlinear optimization problem from \cite{Duran1986}. For demonstrative purposes, the original nonlinear objective has been replaced with a linear objective, and variables $\by$ have been concatenated to $\bx$ for consistency of notation. 

\begin{equation}
    \begin{aligned}
        \text{min} &~ f(\bx) = 10x_1 - 17x_3 -5x_4 + 6x_5 + 8x_6  \\
        \text{s.t.} &~ g_1(\bx) = 0.8\rm{log}(x_2 + 1) + 0.96\rm{log}(x_1 - x_2 + 1) - 0.8x_3 \geq 0, \\ 
        &~ g_2(\bx) = \rm{log}(x_2 + 1) + 1.2\rm{log}(x_1 - x_2 + 1) - x_3 - 2x_6 + 2 \geq 0,  \\
            &~ x_1 - x_2 \geq 0,~~2x_4 - x_2 \geq 0, \\
            &~ 2x_5 - x_1 + x_2 \geq 0,~~1 - x_4 - x_5 \geq 0, \\
            &~ 0 \leq x_1 \leq 2,~0 \leq x_2 \leq 2,~0 \leq x_3 \leq 1, \\
            &~ x_4, x_5, x_6 \in \{0,1\}^3.
    \end{aligned}
    \label{eq:dg}
\end{equation}

 We will focus on the nonlinear inequalities $g_1(\bx) \geq 0$ and $g_2(\bx) \geq 0$ as we implement the method step by step. 

 \subsection{Standard Form Problem.}

 Most global optimization problems are compatible with the standard form in Section~\ref{sec:stdform} by construction. We demonstrate this by partitioning the original problem \eqref{eq:dg} below. 

 \begin{align*}
    \text{min} &~ f(\bx) = 10x_1 - 17x_3 -5x_4 + 6x_5 + 8x_6 & \rm{Objective} & \\    
    \hline
    \text{s.t.} &~ g_1(\bx) = 0.8\rm{log}(x_2 + 1) + 0.96\rm{log}(x_1 - x_2 + 1) - 0.8x_3 \geq 0, & \rm{Nonlinear} & \\ 
    &~ g_2(\bx) = \rm{log}(x_2 + 1) + 1.2\rm{log}(x_1 - x_2 + 1) - x_3 - 2x_6 + 2 \geq 0, & \rm{constraints} & \\
    \hline
        &~ x_1 - x_2 \geq 0,~~2x_4 - x_2 \geq 0, & \rm{Linear}&\\
        &~ 2x_5 - x_1 + x_2 \geq 0,~~1 - x_4 - x_5 \geq 0, & \rm{constraints} & \\
    \hline
        &~ 0 \leq x_1 \leq 2,~0 \leq x_2 \leq 2,~0 \leq x_3 \leq 1, &\rm{Variables} & \\
        &~ x_4, x_5, x_6 \in \{0,1\}^3. & \rm{and~bounds} &
 \end{align*}

 We pass the linear constraints, variables and bounds directly to the \gls{mio} model, and confirm that all variables in nonlinear constraints, in this case $x_1$, $x_2$, $x_3$ and $x_6$, are bounded. Note the presence of binary $x_4$, $x_5$ and $x_6$ in the problem as well. 

 \vspace{0.25in}
 \subsection{Sampling and Evaluation of Nonlinear Constraints.}

\begin{figure}[h!]
     \begin{center}
         \begin{subfigure}{\textwidth}
             \centering
             \includegraphics[width = 0.7\linewidth]{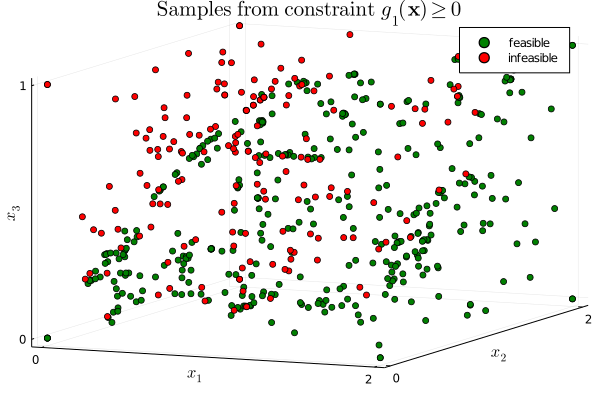}
         \end{subfigure}
     \end{center}
     \caption{The distribution of data for constraint $g_1(\bx) \geq 0$, generated by sampling procedures defined in Section~\ref{sec:sampling}.}
     \label{fig:demosamples}
\end{figure}

 Next we generate samples over the nonlinear constraints using the procedure in Section~\ref{sec:sampling}. Note that $g_1(\bx) \geq 0$ and $g_2(\bx) \geq 0$ have 3 and 4 active variables, so samples are generated in $\R^3$ and $\R^4$ respectively. The resulting samples over $g_1(\bx) \geq 0$ and their feasibilities are shown in Figure~\ref{fig:demosamples}. Note that the samples span the whole $\bx$ hypercube, but that there are certain concentrations of points, thanks to the \gls{knn} sampling procedure, that approximate the constraint boundary. This improves the ability of the approximating \gls{octh} to be both globally and locally accurate.

 \subsection{Decision Tree Training.}

\begin{figure}[h!]
\begin{center}
\begin{subfigure}{0.8\textwidth}
    \centering
    \includegraphics[width = 0.6\linewidth]{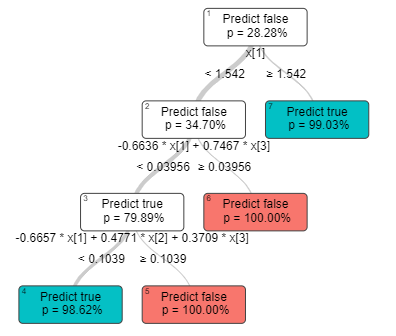}
    \caption{\gls{octh} over constraint $g_1(\bx) \geq 0$ has two feasible and two infeasible leaves, and a depth of 3.}
    \label{fig:tree1}
\end{subfigure} 
\begin{subfigure}{0.8\textwidth}
    \centering
    \includegraphics[width = 0.8\linewidth]{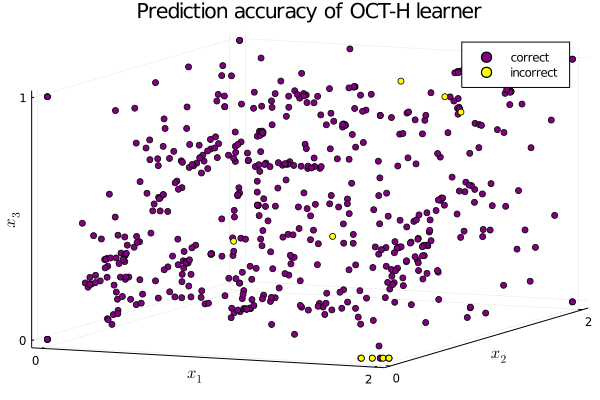}
    \caption{\gls{octh} approximation is 97\% accurate over 554 samples.}
    \label{fig:tree2}
\end{subfigure}
\caption{The approximating OCT-H achieves a high degree of accuracy, capturing both the global and local behavior of the constraint $g_1(\bx) \geq 0$.}
\label{fig:exampletrees}
\end{center}
\end{figure}

We train two \gls{octh}s to classify the feasible space of constraints $g_1(\bx) \geq 0$ and $g_2(\bx) \geq 0$. For demonstrative purposes, the trees were limited to a maximum depth of 3, as opposed to the standard depth of 5 used in \gls{octhagon} as defined in Table~\ref{tab:treeparams}.  The approximating \gls{octh} for $g_1(\bx) \geq 0$ and the accuracy of its predictions are presented in Figure~\ref{fig:exampletrees}. Notably, the \gls{octh} approximator achieves a high degree of accuracy (97\%) throughout $\rm{dom}(\bx)$ with only two feasible leaves. 

\vspace{0.25in}

\subsection{MI Approximation.}

We pose the trees in a \gls{mio}-compatible form. As a bookkeeping note, auxiliary variables are introduced with two indices, the first indicating the constraint index, and the second indicating the numerical index of the leaf of the approximating \gls{octh}. This is consistent with the formulation in Section~\ref{sec:mi_constraints}. 

\begin{figure}[h!]
    \resizebox{\textwidth}{!}{
    \begin{subfigure}{0.5\textwidth}
        \includegraphics[width = \linewidth]{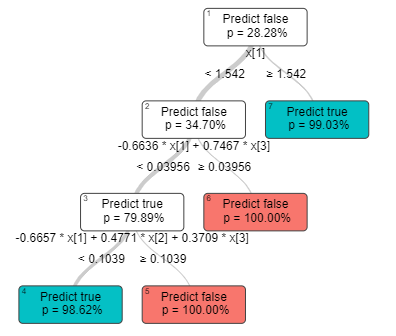}  
    \end{subfigure}
    \begin{subfigure}{0.5\textwidth} 
        \begin{align}
            \label{eq:g1}
            \begin{split}
                [1,0,0] \cdot \by_{1,7} &\geq 1.542 z_{1,7}, \\
                [1,0,0] \cdot \by_{1,4} &\leq 1.542 z_{1,4}, \\
                [-0.6636, 0 , 0.7467] \cdot \by_{1,4} &\leq 0.03956 z_{1,4}, \\
                [-0.6657, 0.4771, 0.3709] \cdot \by_{1,4} &\leq 0.1039 z_{1,4}, \\
        \vspace{0.5in}
                \by_{1,4} + \by_{1,7} = [x_1, x_2, x_3],~&z_{1,4} + z_{1,7} = 1, \\
        \vspace{0.5in}
                [0, 0, 0]z_{1,4} \leq \by_{1,4} &\leq [2,2,1] z_{1,4},\\ 
                [0, 0, 0]z_{1,7} \leq \by_{1,7} &\leq [2,2,1] z_{1,7}, \\ 
                z_{1,4}, z_{1,7} &\in \{0, 1\}^2.
            \end{split}  
        \end{align}
    \end{subfigure}
    }
    \caption{$g_1(\bx) \geq 0$ is approximated via 6 continuous and 2 binary auxiliary variables, and 6 linear constraints.}    
    \label{fig:g1approx}
\end{figure}

\begin{figure}[h!]
    \resizebox{\textwidth}{!}{
    \begin{subfigure}{0.5\textwidth}
        \includegraphics[width = \linewidth]{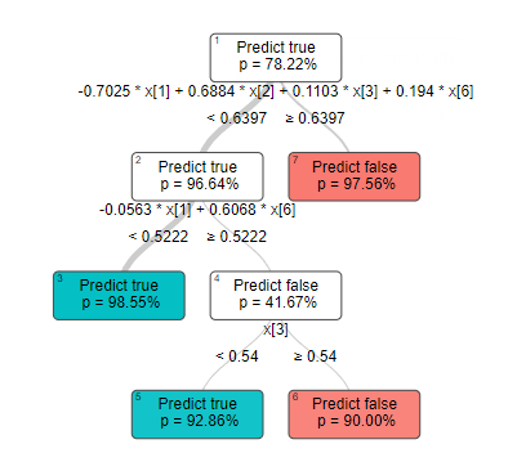}  
    \end{subfigure}
    \begin{subfigure}{0.5\textwidth}
        \begin{align}
            \label{eq:g2}
            \begin{split}
                [-0.7025, 0.6884, 0.1103, 0.194] \cdot \by_{2,3} &\leq 0.6397 z_{2,3}, \\
                [-0.0563, 0, 0, 0.6068] \cdot \by_{2,3} &\leq 0.5222 z_{2,3}, \\
                [-0.7025, 0.6884, 0.1103, 0.194] \cdot \by_{2,5} &\leq 0.6397 z_{2,5}, \\
                [-0.0563, 0, 0, 0.6068] \cdot \by_{2,5} &\geq 0.5222 z_{2,5}, \\
                [0,0, 1, 0] \cdot \by_{2,5} &\leq 0.54 z_{2,5}, \\
                \by_{2,3} + \by_{2,5} = [x_1, x_2, x_3, x_6],~&z_{2,3} + z_{2,5} = 1, \\
                [0, 0, 0, 0]z_{2,3} \leq \by_{2,3} &\leq [2,2,1,1] z_{2,3},\\ 
                [0, 0, 0, 0]z_{2,5} \leq \by_{2,5} &\leq [2,2,1,1] z_{2,5}, \\ 
                z_{2,3}, z_{2,5} &\in \{0, 1\}^2.
            \end{split}  
        \end{align}
    \end{subfigure}
    }
    \caption{$g_2(\bx) \geq 0$ is approximated via 8 continuous and 2 binary auxiliary variables, and 7 linear constraints.}
    \label{fig:g2approx}
\end{figure}

Figure~\ref{fig:g1approx} shows the approximating tree for constraint $g_1(\bx)\geq 0$, as well as its disjunctive representation as defined by \eqref{eq:disjunctive}. Since the constraint has three active variables $[x_1, x_2, x_3]$, and the tree has two feasible leaves with node indices 4 and 7, the disjunctive representation requires the definition of 6 auxiliary continuous variables $\by_{1,4} \in \R^3$ and $\by_{1,7} \in \R^3$, and two binary variables $z_{1,4}$ and $z_{1,7}$. The number of linear constraints required is 6, which is equal to the sum of the depths of each feasible leaf, plus 2 additional constraints defining the disjunctions.

We approximate $g_2(\bx) \geq 0$ in Figure~\ref{fig:g2approx}, with four active variables $[x_1, x_2, x_3, x_6]$, using the same approach. 

\subsection{Solution of MIO Approximation.}

\begin{figure}[h!]
    \begin{subfigure}{0.98\textwidth}
        \begin{center}	
                \begin{tabular}{|c|c|c|c|c|c|c|c|}	
                    \hline	
                    & $x_1$ & $x_2$ & $x_3$ & $x_4$ & $x_5$ & $x_6$ & $f(\bx)$ \\ 	
                    \hline	
                    MIO & 0.375 & 0.375 & 0.379 & 1.0	& 0.0 &	0.0 & -7.685 \\
                    PGD-repaired & 0.699 & 0.699	& 0.530	& 1.0 & 0.0	& 0.0 & -7.021 \\
                    \hline	
                \end{tabular}	
                \caption{The optimal solutions to demonstrative problem, pre- and post-PGD repair.}	
                \label{tab:pgd} 	
        \end{center}
    \end{subfigure}
    \begin{subfigure}{0.98\textwidth}
        \begin{center}
            \includegraphics[width = 0.8\linewidth]{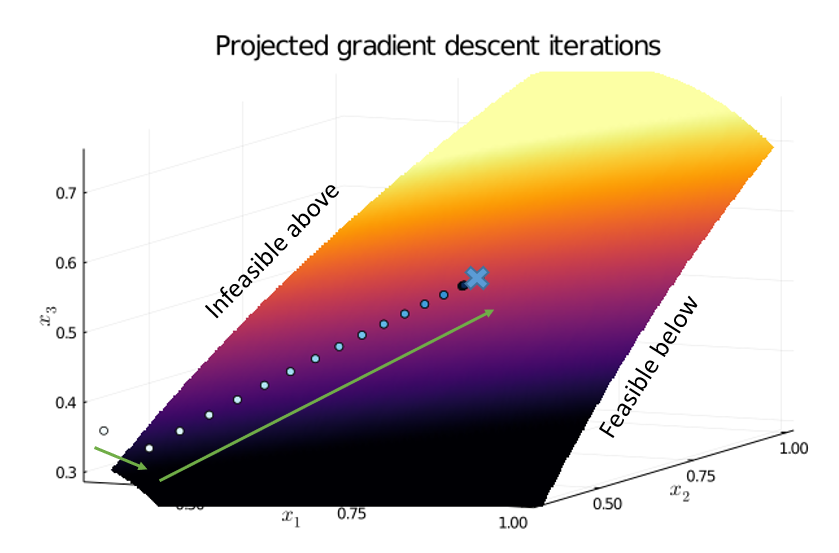}
        \end{center}
        \caption{The progress of the \gls{pgd} method on the demonstrative example, plotted with respect to $x_1$, $x_2$ and $x_3$ on the surface of $g_1(\bx)\geq 0$.}
        \label{fig:pgditerations}
    \end{subfigure}
    \caption{The MIO solution to the demonstrative example is successfully repaired to be feasible and locally optimal by the PGD method.}
    \label{fig:pgdall}
\end{figure}

As described Section~\ref{sec:mio}, once the intractable constraints $g_1(\bx) \geq 0$ and $g_2(\bx) \geq 0$ are replaced with their tractable disjunctive approximations \eqref{eq:g1} and \eqref{eq:g2}, the problem turns into a \gls{milo} that is tractable using commercial solvers. We solve the problem via CPLEX, and obtain a near-feasible, near-optimal solution with the objective value of -7.685 in Table~\ref{tab:pgd}.

\vspace{0.25in}
\subsection{Solution Checking and Repair.}

 We check whether the approximate solution $\bx^{*}$ is feasible to the original optimization problem \eqref{eq:dg} by evaluating the two nonlinear constraints. Since constraint $g_1(\bx)\geq 0$ is violated, we initiate the \gls{pgd} repair procedure from Section~\ref{sec:repair}. To do so, \gls{octhagon} replaces the \gls{mio} approximations of intractable constraints with the auto-differentiated gradients of the constraints at $\bx^*$, and takes a local step to close the feasibility gap while descending along the objective. This is done iteratively, evaluating the objective function and nonlinear constraints at each step, until all constraints are feasible, and the change in the objective value falls below an absolute tolerance ($10^{-4}$). The path of the \gls{pgd} algorithm is shown in Figure~\ref{fig:pgditerations}, on the surface of constraint $g_1(\bx) \geq 0$, which divides the feasible space of $\bx$ in two. Note that this surface is unknown by the method, so it is remarkable that it projects towards it with remarkable accuracy in its first step, and then moves along the surface in a series of descent steps. 
 
 For this problem, the absolute tolerance of $10^{-4}$ was too small to converge definitively, so the \gls{pgd} algorithm terminates at its maximum of 100 iterations, with the optimal objective value of $f(\bx^*) = -7.021$ and the optimal solution in Table~\ref{tab:pgd}.

\section{Computational Experiments on Benchmarks}
\label{sec:benchmarks}

We apply \gls{octhagon} to a number of optimization problems from the literature, and benchmark it against other global optimizers. The software implementation of \gls{octhagon} can be found via the link in Appendix~\ref{app:octhagonimp}. For the full list of optimizers used and their capabilities, please refer to Appendix~\ref{app:solvers}. We lead this section with a caveat. Since our approach is approximate, different random restarts of the solution procedure may yield different optima. However, experience implementing the method suggests that the method is consistent in finding the same optimum in most cases, and that random restarts reliably mitigate issues resulting from finding near-optimal solutions. 

We first apply our method to five small benchmark problems from MINLPLib~[\cite{Bussieck2003}], and compare our results to those of BARON~[\cite{Sahinidis1996}], a popular and effective commercial \gls{minlp} solver. The types and numbers of constraints in the benchmarks are listed in Table~\ref{tab:smallbenchstats}. The results are shown in Table~\ref{tab:smallbenchresults}.

\begin{table}[h]
\begin{center}
        \begin{tabular}{|c|c|c|c|c|c|c|}
            \hline
            \stackedcell[c]{Problem\\Name} &  
            \stackedcell[c]{Continuous\\Variables} & 
            \stackedcell[c]{Integer\\Variables} &
            \stackedcell[c]{Linear\\Constraints} & 
            \stackedcell[c]{Nonlinear\\Inequalities} & 
            \stackedcell[c]{Nonlinear\\Equalities} &
            \stackedcell[c]{Nonlinear\\Objective} \\
            \hline
            \texttt{minlp} & 3 & 1 & 4 & 2 & 0 & Y \\
            \texttt{pool1} & 7 & 0 & 2 & 4 & 0 & N \\
            \texttt{nlp1} & 2 & 0 & 0 & 1 & 0 & N \\
            \texttt{nlp2} & 3 & 0 & 0 & 0 & 3 & N \\
            \texttt{nlp3} & 10 & 0 & 3 & 1 & 3 & Y \\
            \hline
        \end{tabular}
    \caption{The five small nonlinear benchmarks from MINLPLib have a combination of nonlinear inequalities, equalities and objective.}
    \label{tab:smallbenchstats}
\end{center}
\end{table}

\begin{table}[h]
    \begin{center}
        \resizebox{\textwidth}{!}{
            \begin{tabular}{|c|c|c|c|c|c|c|}
                \hline
                Problem name & \multicolumn{2}{|c|}{Objective} & \multicolumn{2}{|c|}{Time (s)} & \multicolumn{2}{|c|}{Solution} \\
                \hline
                & BARON & OCT-HaGOn & BARON & OCT-HaGOn & BARON & OCT-HaGOn \\
                \hline
                \texttt{minlp} & 6.0098 & 6.0098 & 0.120 & 29.9 & [0,1,0,1.3,0,1] & [0,1,0,1.3,0,1] \\
                \texttt{pool1} & 23.0 & 23.0 & 0.082 & 3.90 & 
                \shortstack{[4.0, 3.0, 1.0, 4.0, \\ 0.0  2.12, 0.0]}
                & \shortstack{[4.0, 3.0, 1.0, 4.0, \\
                    0.0  6.63, 0.0]} \\
                \texttt{nlp1} & -6.667 & -6.667 & 0.106 & 0.461 & [6, 0.667] & [6, 0.667] \\
                \texttt{nlp2} & 201.16 & 201.16 & 0.092 & 2.75 & [6.29, 3.82, 201.16] & [6.29, 3.82, 201.16] \\ 
                \texttt{nlp3} & -1161.34 & -1161.34 & 1.265 & 17.7 & [...] & [...] \\
                \hline 
            \end{tabular}
        }
        \caption{Solutions to the small benchmarks using OCT-HaGOn and BARON.}
        \label{tab:smallbenchresults}
    \end{center}
\end{table}

\gls{octhagon} is able to find the global optima for all five small benchmarks, matching the BARON solutions. \gls{octhagon} takes significantly longer to solve the small benchmarks than BARON. This is expected, since these problems have explicit constraints that only contain mathematical primitives BARON supports. Tree training time makes up the vast majority of the solution times for the small benchmarks; the \gls{mio} and \gls{pgd} solution steps are efficient, taking less that 5\% of the total time for each benchmark. Within the context of using optimization in design, where the optimization would be run many times to obtain a number of solutions on the Pareto frontier, \gls{octhagon} is competitive and even faster than BARON, since the \gls{mio} and \gls{pgd} steps are solved in a small fraction of the time it takes for the BARON solver to solve a single instance of each \gls{minlp}.

\begin{table}[h]
    \begin{center}
        \resizebox{\textwidth}{!}{
            \begin{tabular}{|c|c|c|c|c|c|c|}
                \hline
                    \stackedcell[c]{Problem\\name} &
                    \stackedcell[c]{Continuous\\Variables} & 
                    \stackedcell[c]{Integer\\Variables} &
                    \stackedcell[c]{Linear\\Constraints} & 
                    \stackedcell[c]{Nonlinear\\Inequalities} & 
                    \stackedcell[c]{Nonlinear\\Equalities} &
                    \stackedcell[c]{Nonlinear\\Objective} \\
                \hline
                \texttt{himmel16} & 19 & 0 & 1 & 15 & 6 & N \\
                \texttt{kall\_circles\_c6b} & 18 & 0 & 54 & 21 & 1 & N\\
                \texttt{pointpack08} & 17 & 0 & 41 & 28 & 0 & N \\
                \texttt{flay05m} & 23 & 40 & 61 & 5 & 0 & N \\
                \texttt{fo9} & 111 & 72 & 326 & 18 & 0 & N \\
                \texttt{o9\_ar4\_1} & 109 & 72 & 418 & 18 & 0 & N \\
                \hline
            \end{tabular}
        }
        \caption{The six larger benchmarks from MINLPLib. Note that the objective functions are linear in $\bx$, and that nonlinearities are instead embedded in the constraints.}
        \label{tab:largebenchstats}
    \end{center}
\end{table}

\begin{table}[h]
    \begin{center}
        \resizebox{\textwidth}{!}{
            \begin{tabular}{|c|c|c|c|c|c|c|}
                \hline
                Problem name & \multicolumn{3}{|c|}{Objective} & \multicolumn{2}{|c|}{Time (s)} & GO \\
                \hline
                 & GO & OCT-HaGOn & BK & Global & OCT-HaGOn & \\
                 \hline
                \texttt{himmel16} & -0.6798 & $-0.8660^*$ & -0.8660 & 0.055 & 109.575 & CONOPT \\
                \texttt{kall\_circles\_c6b} & 2.8104 & $2.1583^*$ & 1.9736 & 0.355 & 38.503 & IPOPT\\ 
                \texttt{pointpack08} & -0.2574 & -0.2500 & -0.2679 & 13.483 & 91.805 & IPOPT \\
                \texttt{flay05m} & 64.498 & 64.499 & 64.498 & 0.212 & 9.515 & CONOPT \\
                \texttt{fo9} & 23.464 & 23.464 & 23.464 & 959.090 & 29.534 & BARON \\
                \texttt{o9\_ar4\_1} & 236.138 & 236.138 & 236.138 & 2283.281 & 1255.598 & BARON \\
                \hline
            \end{tabular}
        }
        \caption{Solutions to the larger benchmarks using commercial global optimizers (GOs) and OCT-HaGOn, against best known (BK) solutions.}
        \label{tab:largebenchresults}
    \end{center}
\end{table}

We proceed by considering a set of six larger benchmarks from MINLPLib, as shown in Table~\ref{tab:largebenchstats}. We also address the optimization problems using three commercially available solvers, IPOPT, CONOPT and BARON. Given the increased difficulty of these larger benchmarks, we allow \gls{octhagon} to use trees with maximum depth of 8, and generate double the number of samples per constraint compared to the small benchmarks. 

The results are shown in Table~\ref{tab:largebenchresults}, compared with the best known solutions as documented in MINLPLib. \gls{octhagon} finds the best known global optima for 4 out of 6 instances, and high performing solutions otherwise. Some modifications were required for a subset of the problems to be able to apply our method. The problems marked with an asterisk required the following changes to the algorithm:

\begin{itemize}
    \item The \texttt{himmel16} test case contains a number of variables in nonlinear constraints that are unbounded. Using our little knowledge of the problem, we were able to make it compatible with our method by imposing bounds on all variables, $\bx \in [-1, 1]^{19}$. 
    \item The \texttt{kall\_circles\_c6b} example required increasing the step penalty and equality penalty to $10^8$, to damp the \gls{pgd} projection rate in order to avoid a conservative local optimum.
\end{itemize}

While these results are promising in showing that the method can scale to larger problems, they point to some practical considerations. The results show weak correlation between solution time and size of the problems; this is because the number of variables and complexity of nonlinearities in the approximated constraints tend to drive tree training time and thus total solution time. Additionally, as the problem size increases, it is not obvious whether tree training or \gls{mio} steps drive computational time, especially in the presence of integer variables. For the \texttt{himmel16} example, tree training takes 104 seconds of the 110 second total time, whereas for the \texttt{o9\_ar4\_1} benchmark, optimization time dominates, with training only taking 3 seconds out of nearly 21 minutes of total time. And while \texttt{fo9} and \texttt{o9\_ar4\_1} are of similar sizes and have similar constraints (both contain nonlinearities with inverses), they have dramatically different solution times. 

\section{Real World Examples}
\label{sec:realworld}

In addition to the benchmarks, we test our method on two aerospace problems of varying complexity. We first solve a benchmark from the engineering literature, to show that the method can address real-world problems. We then apply \gls{octhagon} to a satellite on-orbit servicing problem that cannot be addressed using other global optimizers. 

\subsection{Speed Reducer Problem.}
\label{sec:speedreducer}

The speed reducer problem is a nonlinear optimization problem posed in~\cite{Golinski1970}. The problem aims to design a gearbox for an aircraft engine, subject to 11 specifications, geometry, structural and manufacturability constraints, in addition to variable bounds over $\bx \in \R^7$. We apply our method to the problem as written in Appendix~\ref{app:sr} in standard form. 

\begin{table}[h!]
\begin{center}
    \resizebox{\textwidth}{!}{
        \begin{tabular}{|r|ccccccc|c|c|c|}
            \hline
            & $x_1$ & $x_2$ & $x_3$ & $x_4$ & $x_5$ & $x_6$ & $x_7$ & Objective & Time (s) & Error\\
            \hline
            BK & 3.5 & 0.7 & 17 & 7.3 & 7.7153 & 3.3503 & 5.2867 & 2994.472 & 476 & $10^{-6}$\\
            \gls{octhagon} & 3.5 & 0.7 & 17 & 7.3 & 7.7153 & 3.3502 & 5.2867 & 2994.355 & 32.6 & 0 \\
            IPOPT & 3.5 & 0.7 & $17.0^*$ & 7.3 & 7.7153 & 3.3502 & 5.2867 & 2994.355 & 4.2 & $10^{-7}$ \\
            \hline
        \end{tabular}
    }
\caption{Both OCT-HaGOn and IPOPT beat the best known (BK) solution of the speed reducer problem. In addition, OCT-HaGOn has 0 error on constraint satisfaction.}
\label{tab:srsol}
\end{center}
\end{table}

In Table~\ref{tab:srsol}, we compare different solutions to the speed reducer problem. Both \gls{octhagon} and IPOPT beat the best known optimum from \cite{Lin2012}. In addition, \gls{octhagon} allows us to achieve all constraints with zero error after 4 iterations of the \gls{pgd} algorithm as shown in Appendix~\ref{app:srpgd}, while the other two methods have small but nonzero error tolerances.

IPOPT was able to solve this particular \gls{nlp} in 4.2 seconds, significantly faster than \gls{octhagon}, which took 32.6 seconds. However, this required a relaxation of the integrality of $x_3$. For this particular problem, this was not concerning since $x_3$ was lower bounded by its optimal value of 17. However, IPOPT cannot in general be used to solve \gls{minlp}s.

On a practical note, we would like to note the different levels of complexity in the \gls{octh} approximations of the underlying nonlinear constraints. Some constraints, while they look quite complex, have low-complexity tree approximators. Consider the following constraint $g_5(\bx) \geq 0$ and its associated \gls{octh} approximator. 

\begin{figure}[p]
        \begin{subfigure}{\linewidth}
            \centering
            \begin{equation*}
                \label{eq:sr_g5}
            g_5(\bx) = 110 x_6^3 - \Big[\Big(745\frac{x_4}{x_2x_3}\Big)^2 + 16.9 \times 10^6\Big]^{0.5} \geq 0
            \end{equation*}
        \end{subfigure} \\
        \begin{subfigure}{\linewidth}
            \centering
            \includegraphics[width=0.55\textwidth]{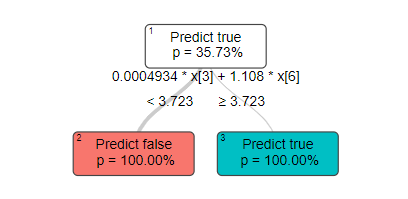}
        \end{subfigure}  
        \caption{The constraint $g_5(\bx) \geq 0$ is accurately approximated by a single separating hyperplane over $\rm{dom}(\bx)$.}
        \label{fig:sr_g5}       
\end{figure}

The \gls{octh} model has a single hyperplane that is able to approximate the function in the relevant $\mathrm{dom}(\bx)$ with perfect accuracy over 613 samples, as shown in Figure~\ref{fig:sr_g5}. Within the bounded $\mathrm{dom}(\bx)$, the nonlinear constraint is thus simplified to a linear constraint. 

\begin{figure}[p]
    \begin{subfigure}{\textwidth}
        \begin{center}
            \begin{align}
                \begin{split}
                    f(\bx) &= 0.7854 x_1 x_2 ^ 2 (3.3333 x_3^2 + 14.9334 x_3 - 43.0934) \\
                    &- 1.5079 x_1 (x_6^2 + x_7^2) + 7.477 (x_6^3 + x_7^3) + 0.7854 (x_4x_6^2 + x_5x_7^2).
                \end{split}
                \label{eq:sr_f}
            \end{align}
        \end{center}
    \end{subfigure}
    \begin{subfigure}{0.9\textwidth}
        \begin{center}
            \includegraphics[width=\linewidth]{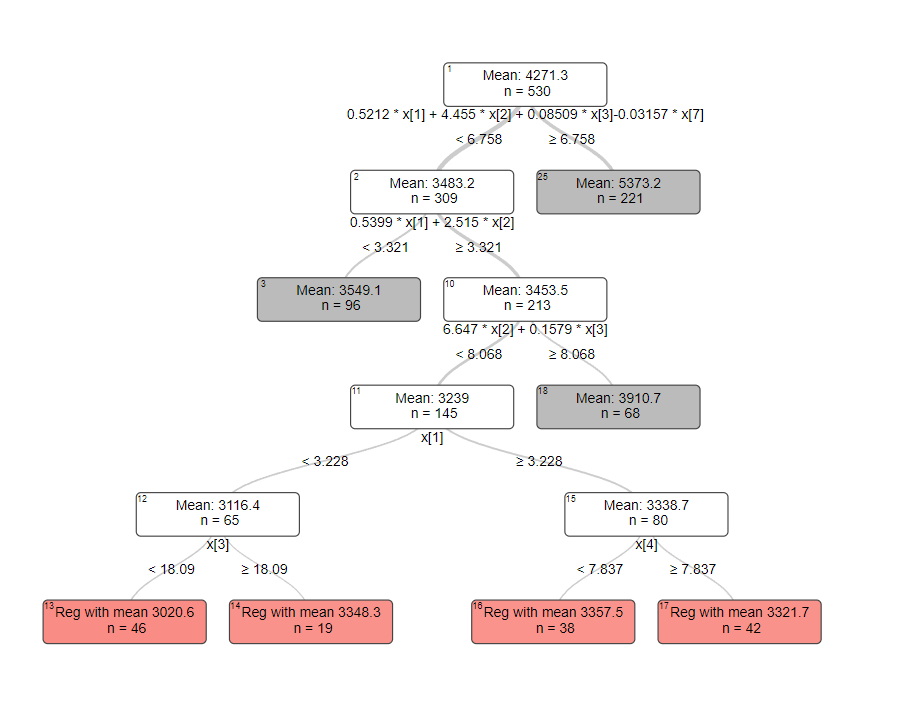}
        \end{center}
    \end{subfigure}
    \caption{The objective function $f(\bx)$ is approximated via an ORT-H with 19 leaves (4 leaves shown) and $1 - \mathrm{R}^2$ error of $1.4\times10^{-5}$.}
    \label{fig:sr_f}
\end{figure} 

However, not all constraints are straightforward to represent via unions of polyhedra. Consider the objective function, which is a 5th order polynomial \eqref{eq:sr_f}. In this particular case, the objective is represented by an \gls{orth} with 19 leaves, each defining a unique feasible polyhedron over $\bx$. A truncated version of the tree, with four leaves visible, is shown in Figure~\ref{fig:sr_f}. The $1-\rm{R}^2$ error of the approximation is $1.4\times10^{-5}$ over 532 samples.

\subsection{Satellite OOS Problem.}
\label{sec:oos}

We test our method on the previously-unsolved optimization problem of satellite \gls{oos} scheduling. Satellite \gls{oos} is a future technology that seeks to improve the lifetime of existing and next-generation satellites by allowing autonomous servicer spacecraft to perform repairs or refuels in orbit [\cite{Luu2020}]. \gls{oos} is a difficult scheduling problem that acts on a highly nonlinear dynamical system. It is a good problem to address via our method since, in its full \gls{minlp} form, the problem is a nonconvex combinatorial optimization problem with nonlinear equality constraints. In addition, due to the 11 orders of magnitude difference in the ranges of decision variables, it is numerically challenging. Before this paper, it was addressed only via enumeration [\cite{Luu2020}]. Please refer to Appendix~\ref{app:oos} for more details on the full list of constraints; a succinct summary of the problem follows. 

The dynamical problem is the orbital mechanics of moving a servicer satellite between client satellites in the same orbital plane. Orbital transfers involve using on-board thrusters to get the servicer into a different orbital altitude than the client satellite, called the phasing orbit, in order to reduce the true anomaly (angular phase difference in radians) between the servicer and the client. The servicer then propels itself back onto the client's orbit to meet the client satellite at the right time and position in space, while obeying conservation of energy, momentum and mass. The scheduling problem involves both choosing the optimal order in which to serve each client satellite (discrete decisions), as well as choosing the optimal phasing orbits (continuous decisions). 

In this section, we consider a simple example of \gls{oos}. We schedule a single servicer satellite to refuel 7 client satellites in orbit, traveling between clients using on-board propulsors. Each client requires different amounts of fuel, and we constrain the servicer to fulfill its mission in 0.35 years, with the objective being to minimize the wet mass (the dry mass and fuel) of the servicer.  The problem parameters are in Table~\ref{tab:oosparams}.
\begin{table}[h!]
\begin{center}
        \centering
        \begin{tabular}{|c|c|c|}
            \hline
            Parameter & Value & Units \\
            \hline
            Servicer dry mass & 500 & kg \\
            Propulsor specific impulse & 230 & (Ns)/kg \\
            Number of client satellites & 7 & - \\
            Client satellite altitude & 780 & km \\
            Servicer satellite altitude range & [760,800] & km \\
            Maximum service time & 0.35 & years \\
            \hline
        \end{tabular}
        \caption{OOS problem parameters.}
        \label{tab:oosparams}
\end{center}
\end{table}

\begin{center}
\begin{figure}[!h]
        \centering
        \includegraphics[width = 0.65\linewidth]{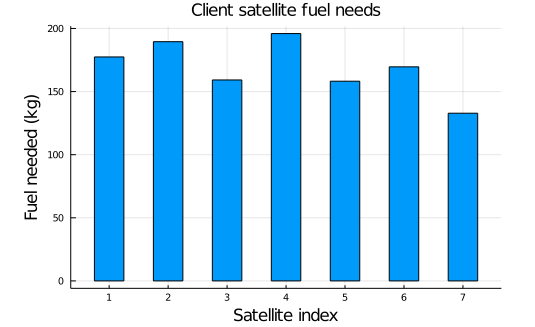}
        \caption{Client satellites require different amounts of fuel, which affects the optimal schedule for servicing.}
        \label{fig:clientSatFuel}
\end{figure}
\end{center}

The fuel requirements shown in the Figure~\ref{fig:clientSatFuel} were randomly generated and reflect a possible distribution of fuel needs for client satellites that are part of the same constellation and were launched concurrently at a previous point in time.  

In addressing the \gls{oos} problem, we make the following realistic simplifying assumptions, although our method does extend to more general cases.  The servicer satellite is delivered by an external rocket to the first client, and uses its own propulsor to use Hohmann transfers between the subsequent client satellites. Thrusting and refueling steps take a negligible amount of time relative to maneuver steps. All client satellites are in the same orbital plane, at the same altitude, and are evenly spaced around the orbit. 

The initial problem of servicing $n_s = 7$ clients has 141 variables, of which $n_s^2 = 49$ binary variables denote the servicing order. The continuous decision variables in nonlinear constraints are bounded from above and below to be compatible with \gls{octhagon} as defined in Section~\ref{sec:stdform}. There are 41 linear constraints in the model representing a subset of the system dynamics.  On top of the linear constraints, we have $10(n_s-1) = 60$ nonlinear constraints, all of them equalities. The constraints are presented in detail in Appendix~\ref{app:oos}.

We solve the problem in two ways. First we solve it via \gls{octhagon}. Since we know the constraints of this problem explicitly, we use the \gls{orth} approximation method as described in Section~\ref{sec:nlobjectives}, separating nonlinearities from affine components of constraints for improved accuracy, and training a tree for each set of recurrent constraints. The resulting \gls{mio} problem has 999 continuous and 349 binary variables, and 3650 linear inequalities and 286 linear equalities. 

Other global optimizers such as CONOPT, IPOPT and BARON cannot be used as benchmarks for \gls{octhagon} on this particular problem. Since \gls{oos} is a mixed-integer problem, gradient-based optimizers such as CONOPT or IPOPT are rendered ineffective, and BARON does not support the nonlinearities present in orbital dynamics. Instead, we successfully discretize out a subset of the nonlinearities in constraints by restricting the possible transfer orbits into 1 km bins. This reduces the complexity of the \gls{oos} problem to a \gls{mi}-bilinear problem, which we are able to solve via Gurobi's \gls{mi}-bilinear optimizer~[\cite{gurobi2021}]. The \gls{mi}-bilinear representation has 394 variables, of which 289 variables are binary. 36 of the 60 nonlinear constraints are turned into bilinear equalities, while the rest are transformed into linear constraints. The solution of the discretized problem is globally optimal, but guaranteed to be worse than the global optimum of the full \gls{minlp} formulation, since a discrete set of orbit altitudes is more restrictive than a continuous set. However, the solution is granular enough to be a good benchmark for \gls{octhagon}.

\begin{table}[h!]
\begin{center}
        \begin{tabular}{|c|c|c|c|c|c|c|c|}
            \hline
            Metric & \multicolumn{7}{|c|}{Values} \\
            \hline
            \multicolumn{8}{|c|}{\gls{octhagon} solution} \\
            \hline
            Wet mass (kg) & \multicolumn{7}{|c|}{1725.9} \\
            Total maneuver time (years) & \multicolumn{7}{|c|}{0.350} \\

            \hline
            Satellite order & 4 & 3 & 2 & 1 & 7 & 6 & 5 \\
            Refuel mass (kg) & 196.0 & 159.2 & 189.5 & 177.4 & 132.9 & 169.6 & 158.2 \\
            Transfer orbit altitude (km) & & 765.8 & 765.8 & 765.8 & 765.8 & 765.8 & 767.6 \\            
            Maneuver fuel (kg) & & 9.60 & 8.74 & 7.73 & 6.79 & 6.08 & 4.17 \\
            Maneuver time (days) & & 20.7 & 20.7 & 20.7 & 20.7 & 20.7 & 24.1 \\
            Orbital revolutions & & 297.0 & 297.0 & 297.0 & 297.0 & 297.0 & 345.3 \\
            \hline
            \multicolumn{8}{|c|}{Discretized \gls{mi}-bilinear solution} \\
            \hline
            Wet mass (kg) & \multicolumn{7}{|c|}{1724.4} \\
            Total maneuver time (years) & \multicolumn{7}{|c|}{0.350} \\
            \hline
            Satellite order & 4 & 3 & 2 & 1 & 7 & 6 & 5 \\
            Refuel mass (kg) & 196.0 & 159.2 & 189.5 & 177.4 & 132.9 & 169.6 & 158.2 \\
            Transfer orbit altitude (km) & & 768.0 & 768.0 & 766.0 & 765.0 & 765.0 & 762.0 \\
            Maneuver fuel (kg) & & 8.46 & 7.53 & 7.51 & 6.77 & 5.80 & 5.51 \\
            Maneuver time (days) & & 24.9 & 24.9 & 21.4 & 19.9 & 19.9 & 16.6 \\
            Orbital revolutions & & 357.1 & 357.1 & 306.1 & 285.7 & 285.7 & 238.1 \\
            \hline
        \end{tabular}
    \caption{The discretized and OCT-HaGOn formulations come up with the same optimal satellite schedule, although the discretized solution performs $0.1\%$ better.}
    \label{tab:oosresults}
\end{center}
\end{table}

\begin{figure}
    \begin{center}
        \begin{subfigure}{0.8\linewidth}
            \includegraphics[width=\textwidth]{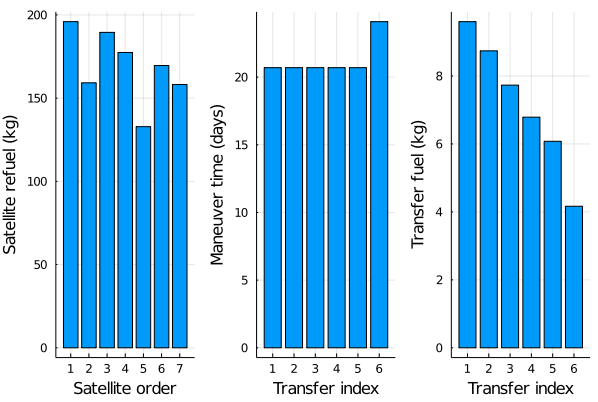}
            \caption{The \gls{octhagon} solution.}
            \label{fig:oosoct}
        \end{subfigure} \\
        \begin{subfigure}{0.8\linewidth}
            \includegraphics[width=\textwidth]{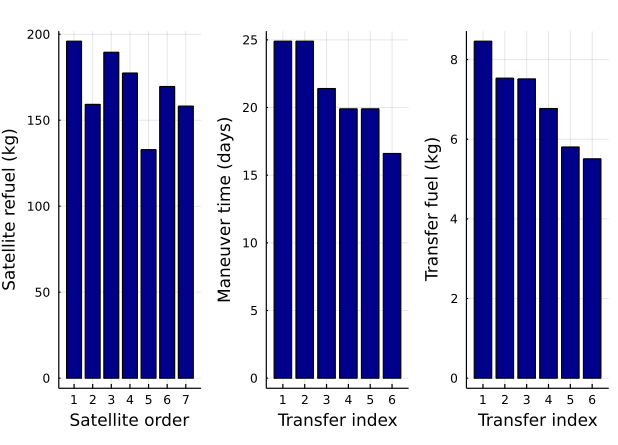}
            \caption{The \gls{mi}-bilinear solution.}
            \label{fig:oosmi}
        \end{subfigure} 
    \end{center}
    \caption{While it captures the orbital dynamics well, OCT-HaGOn is not able to schedule the phasing orbits as well as the MI-bilinear formulation.}
    \label{fig:octcompare}
\end{figure}

The results are presented in Table~\ref{tab:oosresults}, and shown graphically in Figure~\ref{fig:octcompare}. Firstly, we look for two important effects, demonstrated well by the \gls{mi}-bilinear solution and easily seen in Figure~\ref{fig:octcompare}. The first is that it is best to refuel satellites with the largest refuel requirements first, since a lighter servicer requires less fuel to transfer between subsequent clients. The second is that it is better to spend more time transferring in the beginning of the mission than the end, since transfers spend less fuel when the servicer is lighter. This is exhibited by a general downward trend in both maneuver times and fuel costs in the \gls{mi}-bilinear solution. 

While \gls{octhagon} properly captures the optimal satellite schedule, it isn't able to find the optimal set of phasing orbits. This is easily seen by observing the flat profile of maneuver times in the \gls{octhagon} solution in Figure~\ref{fig:oosoct}, which is suboptimal (by $<0.1\%$ total fuel) to the decreasing profile seen in the discretized solution in Figure~\ref{fig:oosmi}. In addition, due to the presence of many nonlinear equalities, the \gls{pgd} method was not able to reduce the infeasibility and optimality gaps, getting stuck in a local optimum. With a maximum tree depth of 6, the solution has a maximum relative error of $3.5 \times 10^{-3}$ and a mean relative error of $2.5 \times 10^{-4}$ on all nonlinear constraints. While this is sufficiently accurate for conceptual design purposes, greater accuracy and a more robust repair procedure are desired. 

In terms of solution time, \gls{octhagon} took 14.2 seconds when solved using a personal computer with an 8-core Intel i7 processor. That includes all sampling, evaluation, training and optimization steps. In comparison, the \gls{mi}-bilinear solution took 17.7 seconds, just for the optimization step. This is in addition to the two days spent by an experienced engineer, reformulating the problem to be compatible with existing efficient optimization formulations.

Despite the suboptimal solution of \gls{octhagon} to the \gls{oos} problem, we argue that it is a strong demonstration of the capabilities and promise of the method, especially considering the problem complexity. Notably, \gls{octhagon} successfully finds the optimal satellite servicing schedule, which is arguably the most important decision in the problem. This is despite the fact that the problem is ill-conditioned, with 11 orders of magnitude difference in decision variable values, and has 60 nonlinear equality constraints coupling a majority of the decision variables. In addition, discretized reformulations of such complex global optimization problems may not exist in general. Even if they do, they may be intractable due to the combinatorial nature of such reformulations. To the best of our knowledge, this makes \gls{octhagon} the only global optimization tool in the literature that can address this problem directly. 

\section{Discussion}
\label{sec:godiscussion}

In this section, we discuss the results and limitations, and propose areas for future work. 

\subsection{Limitations.}

The proposed method shows promise in solving a variety of global optimization problems, but it is a work in progress. Here we detail some of the limitations of the method as implemented in this paper; we list these in order from the most to the least significant, in the author's opinion. 

While the \gls{oos} example demonstrates that the method can address problems with a high degree of nonlinear coupling between decision variables, individual nonlinear constraints involving a large number of active variables will pose challenges in both the \gls{octh} training time, as well as the accuracy of the tree approximations. Tree accuracy directly affects the quality of the approximate optima. Separability, as described in Section~\ref{sec:mi_constraints}, can partially mitigate this problem, by allowing many nonlinear constraints to be decomposed into linear components and better approximated via a series of \gls{orth}s. 

In addition, we have yet to rigorously test how solution time and quality scale with the number of variables and nonlinear constraints, and the sparsity of the nonlinear constraints. Given that the performance of \gls{octhagon} is formulation-dependent, there is much to be gained, both in terms of solution time and quality, through formulations that premeditate where \gls{octh} approximations need to be used, and use them judiciously. We expect \gls{octhagon} to be particularly effective when a majority of the constraints in the optimization problem are linear or convex and therefore efficient, and the constraint learning approach is implemented on the otherwise intractable constraints, with reasonably tight decision variable bounds. 

As noted in Section~\ref{sec:benchmarks}, the proposed method has no guarantees of global optimality since it is approximate. Thus, different iterations of the method generate high-performing solutions that are locally optimal, but do not have guarantees of global optimality such as those provided by BARON. In addition, while the method is agnostic about whether constraints are explicit or inexplicit, the method has so far been tested on explicit constraints only. This is because of the inavailability of numerical benchmarks with black box functions, due to their incompatibility with other existing global optimizers. 

An implicit assumption of the method is that the intractable constraints are quick to evaluate; if this assumption is not true, then the implementation may need to change to accommodate computational requirements. Additionally, the \gls{pgd} method requires that the constraint functions are auto-differentiable. While this is a modest assumption, it is possible that constraint evaluations do not allow for \gls{ad}. This could be overcome by finding gradients approximately, e.g. via finite differencing, but this is not currently implemented. 

With these limitations outlined, we continue by proposing future work to improve the method. 

\subsection{Decision Tree Training.}

The majority of the solution time of \gls{octhagon} is taken by the tree training step. While the computational cost of training is linear with the number of constraints, the results on benchmarks in Section~\ref{sec:benchmarks} show that training time can vary dramatically depending on the complexity of the underlying constraints. In this section, we discuss several ways to manage computational time. 

The first potential source of training time reduction comes from tuning the base tree parameters described in Table~\ref{tab:treeparams} and implemented in \gls{iai}. To do so, we can reduce the complexity of the trees, by reducing the maximum depth and increasing the minbucket parameters. Otherwise we can modify the number of random restarts in tree training. Since the local search method used in generating \gls{octh}s and \gls{orth}s is locally optimal, we can reduce training time by changing the number of random restarts of candidate trees, as well as the number of random hyperplane restarts. However, both methods have a clear negative tradeoff with respect to the accuracy of the \gls{octh} approximations. In general, we find that using 10 random tree restarts and 5 hyperplane restarts, as described by the base tree parameters in Table~\ref{tab:treeparams}, we are able to generate trees that are sufficiently accurate for decision making while being efficient enough to use in a real-time optimization setting. 

A potentially large source of training time reduction is from recognizing the common form of constraints in a problem. If a nonlinear constraint $g(\bu_i) \geq 0$ is repeated $k$ times with different variables $\bu_i \subset \bx,~i \in [k]$, the constraints can be approximated jointly. Specifically, we can train a single \gls{octh} to approximate the constraint over the domain $\cup_{i=1}^k \rm{dom}(\bu_i)$. We then express the $k$ constraints as $k$ repetitions of the disjunctive representation of the tree with different variables $\bu_i$. In this paper, many benchmarks in Section~\ref{sec:benchmarks} exhibit this kind of repeating behavior, but we treat the constraints as black boxes and do not take advantage of potential speed-ups. For the \gls{oos} problem however, we use our knowledge of the constraints to train the trees jointly. 

There is also potential in exploiting the noise-free nature of data over explicit constraints to speed up the training process. Currently, training time scales exponentially in the number of features of the data (i.e. number of variables in each constraint), making tree approximations of constraints dense in $\bx$ slow. One could speed up the training process by trying a greedy approach, building trees with hyperplanes in a locally optimal manner similar to CART [\cite{Breiman1984}], instead of a globally optimal manner via local search heuristics [\cite{Bertsimas2019}]. Another approach could devise specific local search heuristics that sample constraints and train trees in a dynamic manner, in order to speed up training and also reduce the approximation error.

There are also improvements that could be made considering computing architecture. Since individual constraints are learned separately, the training process could be done in parallel, making the best of use of available computational resources. The trees can be efficiently stored once trained, allowing the same trees to be used in different instances of the same optimization problem. This avoids the need to retrain trees, and also avoids having to store the samples required to train the trees, saving on memory. We have developed such methods for development purposes.

\vspace{0.25in}
\subsection{Complexity of the MIO Approximation.}

As aforementioned, the complexity of solving the \gls{mio} approximations of global optimization problems is modest, since the scale of the \gls{mio} is small compared to the abilities of commercial solvers such as Gurobi or CPLEX. However, it is important to note how the complexity of the \gls{mio} can scale depending on the number of nonlinear constraints and the depth of the approximating trees. 

We first consider the number of auxiliary variables required to pose the \gls{mio} approximation. The number of variables used to approximate a nonlinear constraint is a linear function of the number of disjunctive polyhedra describing the feasible space of $\bx$, as well as the number of decision variables in the constraint. More explicitly, the total number of binary variables required to approximate the problem is linear with respect to the number of leaves in the decision trees, and equivalent to
\begin{equation*}
    |L_f| + \sum_{i \in I} |L_{i,1}| + \sum_{j \in J} |L_{j}|,
\end{equation*}
where $L_f$, $L_{i,1}$ and $L_{j}$ are the set of feasible leaves in the objective-, inequality- and equality-approximating trees respectively. In addition, we introduce a number of continuous auxiliary variables. The number of auxiliary variables is equivalent to:
\begin{equation*}
    1 + |L_f|(p_f + 1) + \sum_{i \in I} \Big(|L_{i,1}| p_i\Big) + \sum_{j \in J} \Big(|L_{j}| p_j \Big),
\end{equation*}    
where $p_i$ is the number of variables in the $i$th constraint. The maximum number of leaves of a tree is $2^d$, so in the worst case, the number of auxiliary binary variables in the problem is $\mathcal{O}(2^d(1+|I|+|J|))$, and the number of auxiliary continuous variables is $\mathcal{O}(2^d(1+|I|+|J|)\rm{dim}(\bx))$, equivalent to the number of binary variables augmented by the dimension of $\bx$. In practice however, this worst case is not seen, as the trees are pruned during the training process, and approximated intractable constraints are sparse in $\bx$. 

The number of disjunctive constraints is more complicated, since the trees are not guaranteed to be of uniform depth, and we do not know a priori the fraction of feasible leaves for a classification tree. However, if we assume that each tree has a depth $d_i$, we get the following worst case number of disjunctive constraints, not including the univariate bounding constraints for the continuous auxiliary variables:
\begin{equation*}
    \Big(2^{d_f} \times (d_f + 1)\Big) + 3 + \sum_{i \in I} \Big((2^{d_i}-1) \times d_i \Big) + \sum_{j \in J} \Big((2^{d_j}-1) \times d_j \Big) + 2|I| + 4|J|.
\end{equation*}
The above implies that the number of disjunctive constraints in the \gls{mio} is $\mathcal{O}(2^dd(1+|I|+|J|))$, where $d$ is the maximum depth of all approximating trees. This shows the super-exponential impact of tree depth on \gls{mio} complexity, where the need for greater accuracy may result in large computational cost. However, for the small to medium scale instances we have considered in this paper, this is an acceptable tradeoff.

Additionally, the number of variables grows linearly with number of constraints, which could result in the solution time of \gls{octhagon} being exponential in the worst-case. Unlike linear or convex optimization problems, where the average solution time can be sublinear with the number of constraints, \gls{octhagon} is expected to have on average super-linear solution time with respect to number of constraints due to the combinatorial nature of the approximations. We have yet to observe problems that exhibit such exponential-time behavior, likely because of the sparsity of the approximating constraints, and also due to the locally-idealness of the formulation. However, tree complexity needs to be investigated as \gls{octh} approximations are applied to large scale problems.  

\vspace{0.125in}
\subsection{Extending to MI-Convex Formulations.}

The \gls{octhagon} approach allows us to generate efficient \gls{mio} representations of nonlinear constraints that are not efficiently optimizable, i.e. not linear or convex. It opens up the possibility to include these approximations in more general \gls{mi}-convex formulations, where the efficient convex nonlinear constraints are preserved, either via direct insertion or via outer approximation, while the intractable constraints are approximated via \gls{octh}s. This will significantly improve both the speed and accuracy of our method. 

\vspace{0.125in}
\subsection{Comparing the Big-M Free and Big-M Disjunctive Formulations.}

While \gls{octhagon} implements the locally ideal, big-M free disjunctive representations of decision trees as described in Section~\ref{sec:mi_constraints}, it is possible that a big-M representation is faster to solve via commercial solvers, due to the large number of auxiliary variables added in the big-M free approach. It remains to be tested whether it is more efficient to solve the locally ideal but much larger \gls{milo} resulting from the big-M free approach, or whether it is more efficient to solve the smaller but non-ideal \gls{milo} resulting from the big-M approach. While we have no definitive proof of the relative performance of the two approaches, the author's intuition would point towards a tradeoff based on problem size; it is likely that a big-M approach will outperform the big-M free approach when addressing larger global optimization problems with more nonlinear constraints. 

\vspace{0.125in}
\subsection{Improved Random Restarts.}

As aforementioned, since the constraint learning approach is approximate, random restarts may be required gain confidence in the quality of the locally optimal solutions. Currently, random restarts for \gls{octhagon} involve retraining trees over all nonlinear constraints, and replacing them simultaneously. A better method would be to train an ensemble of trees on each constraint, and permute the tree approximations to generate a set of \gls{mi} approximations of the problem. The solution of each permutation would provide a near-optimal seed for a new \gls{pgd} sequence. This would reduce the computational burden of random restarts and result in higher-performing populations of solutions, giving increased confidence in the method. 

\vspace{0.125in}
\subsection{Optimization Over Data-Driven Constraints.}

There are global optimization contexts where constraints are informed by data, without having access to the underlying models. Some examples are simulation data in the design of engineered systems, outcomes of past experiments, or anthropogenic data such as clinical data and consumer preferences. In theory, \gls{octhagon} is able to learn constraints from arbitrary data and integrate these models in an optimization setting. However, we have yet to perform experiments to confirm the efficacy of \gls{octhagon} in real-world decision making using data-driven constraints. Such an embedding of data into optimization via constraint learning has important implications for a variety of fields, such as healthcare and operations research. 

\vspace{0.125in}
\subsection{Integration of Other MIO-Compatible ML Models.}

While this paper focuses on the use of \gls{octh}s and \gls{orth}s for constraint learning, there are other ML models that have optimization-compatible representations. \cite{Maragno2021} explore the possibility of using linear models, decision trees and their variants, and multi-layer perceptrons to learn constraints and objectives from data. \gls{octhagon} could easily be extended to accommodate such other \gls{mio}-representable \gls{ml} models.

\section{Conclusion}
\label{sec:goconclusion}

In this paper, we have proposed an intuitive new method for solving global optimization problems leveraging interpretable \gls{ml} and efficient \gls{mio}. Our method approximates explicit and inexplicit nonlinear constraints in global optimization problems using \gls{octh}s and \gls{orth}s, using the natural disjunctive representation of decision trees. We demonstrate, both theoretically and practically, that the disjunctive \gls{mio} approximations are efficiently solvable using modern solvers, and result in near-optimal and near-feasible solutions to global optimization problems. We then improve our solutions using gradient-based methods to obtain feasible and high-performing solutions. We demonstrate that our global optimizer \gls{octhagon} is competitive with other state-of-the-art methods in solving a number of benchmark and real-world problems. The Julia implementation of \gls{octhagon} as described in this paper is available via the link in Appendix~\ref{app:octhagonimp}.

The method we present is more than a new tool in the global optimization literature. Tree-based optimization stands out among existing global optimization tools because it can handle constraints that are explicit and inexplicit, and even learn constraints from arbitrary data. To the author's best knowledge, it is the most general global optimization method in the literature, since it has no requirements on the mathematical primitives of constraints or variables. Our method only requires a bounded decision variable domain over the nonlinear constraints. This has important implications to a number of fields that can benefit from optimization, but have yet to do so due to lack of efficient mathematical formulations.

% % UNCOMMENT HERE...
% \section{Acknowledgements}

% We thank Michael Luu for his help in the formulation of the satellite \gls{oos} problem in Section~\ref{sec:oos}, based on his previous work~[\cite{Luu2020}]. 

% This research was supported in part by the Defense Advanced Research Projects Agency (DARPA), under contract HR0011-17-2-0028. The views, opinions and/or findings expressed are those of the author and should not be interpreted as representing the official views or policies of the Department of Defense or the U.S. Government.

% \bibliographystyle{plain}
% \bibliography{main}

% \section{Appendices}
% \input{appendix.tex}

% \end{document}
% % TO HERE! YOU ARE DONE. 

% Appendix here
\begin{APPENDICES}
\section{Appendices}
\subsection{OCT-HaGOn Implementation.}
\label{app:octhagonimp}

\gls{octhagon} is implemented in Julia 1.5.4 available for use at \\
\url{https://1ozturkbe.github.io/research}. The current implementation requires an academic license for Interpretable AI~[\cite{InterpretableAI}], but a lightweight version without Interpretable AI is also in development. While CPLEX is \gls{octhagon}'s default solver, it also supports other \gls{mio} solvers that are compatible with JuMP.jl version 0.21.5~[\cite{DunningHuchetteLubin2017}]. 

\subsection{Optimizers.}
\label{app:solvers}

We use a variety of commercially available and free solvers to address different types of optimization problems. This appendix provides a quick overview of the different optimization tools, the versions used and their capabilities as of writing, as well as their specific applications to different problems.

\begin{itemize}
    \item \textbf{CPLEX v20.1.0.0}: CPLEX, short for ILOG CPLEX Optimization Studio, is a mixed-integer convex optimizer. It is the default solver of \gls{octhagon}, since CPLEX is available for free to solve problems with up to 1000 variables and constraints. In addition, academics can get an unlimited, no-cost academic license. CPLEX is used within \gls{octhagon} to solve the tree-based \gls{mi} approximations of global optimization problems, as well as the \gls{mi}-quadratic optimizations required for the \gls{pgd} iterations. CPLEX is also used in the machinery of BARON, another global optimizer; see below.  
    \item \textbf{Gurobi v9.1.1}: Gurobi is a mixed-integer convex optimizer [\cite{gurobi2021}]. Gurobi is available at no cost via an academic license. Due to its ability to address mixed-integer bilinear optimization problems, Gurobi was used to solve the discretization of the \gls{oos} problem in Section~\ref{sec:oos}, as a benchmark for \gls{octhagon}.
    \item \textbf{CONOPT v3.10:} CONOPT is a gradient-based nonlinear optimizer [\cite{ArneS.Drud1994}]. It was used to solve two large benchmarks in Section~\ref{sec:benchmarks}, via a one-year demo license obtained through the General Algebraic Modeling System (GAMS) interface.  
    \item \textbf{IPOPT v3.13.4:} IPOPT is a freely available interior point optimizer for \gls{nlp}s [\cite{Wachter2006}]. It was used in Section~\ref{sec:benchmarks} to solve two large benchmarks, and in Section~\ref{sec:speedreducer} to address the speed reducer problem. 
    \item \textbf{BARON v21.1.13:} BARON is a commercially available \gls{minlp} solver that accepts a subset of nonlinear primitives~[\cite{sahinidis:baron:21.1.13}]. BARON uses CPLEX as its back-end \gls{mio} solver for its branch-and-reduce solution approach. We purchased a BARON license to be able to solve 5 small benchmarks and 2 large benchmarks in Section~\ref{sec:benchmarks}. 
\end{itemize}

\subsection{Speed Reducer Problem.}
\label{app:sr}

We detail the constraints in the speed reducer problem addressed in Section~\ref{sec:speedreducer}. Note that it has been transcribed from~\cite{Ray2003} into standard form as defined in Section~\ref{sec:stdform}.
\begin{align*}
    \begin{split}
            \underset{\bx}{\text{min}} &~~ 0.7854x_1x_2^2(3.3333 x_3^2 + 14.9334x_3 - 43.0934)  \\
            &~~~~~~ - 1.5079x_1(x_6^2 + x_7^2) + 7.477(x_6^3 + x_7^3)\\
            \text{s.t.} &~~ -27 + x_1  x_2^2  x_3 \geq 0,~-397.5 + x_1  x_2^2  x_3^2 \geq 0, \\
            &~~ -1.93 + \frac{x_2  x_6^4  x_3}{x_4^3} \geq 0,~-1.93 + \frac{x_2  x_7^4  x_3}{x_5^3} \geq 0, \\
            &~~ 110.0x_6^3 - \Bigg(\Big(\frac{745x_4}{x_2x_3}\Big)^2 +
            16.9 \times 10^6\Bigg)^{0.5} \geq 0, \\
            &~~ 85.0x_7^3 - \Bigg(\Big(\frac{745x_5}{x_2x_3}\Big)^2 +
            157.5 \times 10^6\Bigg)^{0.5} \geq 0, \\
            &~~ 40 - x_2x_3 \geq 0,~x_1 - 5x_2 \geq 0,~12x_2 - x_1 \geq 0, \\
            &~~ x_4 - 1.5x_6 - 1.9 \geq 0,~x_5 - 1.1x_7 - 1.9 \geq 0, \\
            &~~ \bx \geq [2.6, 0.7, 17, 7.3, 7.3, 2.9, 5], \\
            &~~ \bx \leq [3.6, 0.8, 28, 8.3, 8.3, 3.9, 5.5], \\
            &~~ x_3 \in \mathbb{Z}.
    \end{split}
\end{align*}

\subsubsection{Speed Reducer PGD Iterations.}
\label{app:srpgd}

The speed reducer problem is converged to a feasible and locally optimal solution from the \gls{mio} solution in 4 \gls{pgd} steps. The decision variable and objective values at each iteration are given in Table~\ref{tab:srpgd}.
\begin{table}
    \begin{center}
        \begin{tabular}{|c|cccccccc|}
            \hline
            Iteration & $x_1$ & $x_2$ & $x_3$ & $x_4$ & $x_5$ & $x_6$ & $x_7$ & Objective\\
            \hline
            1 & 3.5 & 0.7 & 17.0 & 7.3 & 7.71590 & 3.35011 & 5.28718 & 3018.809 \\
            2 & 3.5 & 0.7 & 17.0 & 7.3 & 7.71590 & 3.35011 & 5.28718 & 2994.674 \\
            3 & 3.5 & 0.7 & 17.0 & 7.3 & 7.71590 & 3.35021 & 5.28718 & 2994.700 \\
            4 & 3.5 & 0.7 & 17.0 & 7.3 & 7.71532 & 3.35021 & 5.28665 & 2994.355 \\
            5 & 3.5 & 0.7 & 17.0 & 7.3 & 7.71532 & 3.35021 & 5.28665 & 2994.355 \\
            \hline
        \end{tabular}
    \end{center}
    \caption{Speed reducer PGD iterations.}
    \label{tab:srpgd}
\end{table}

\subsection{Satellite Dynamics Optimization.}
\label{app:oos}

The satellite \gls{oos} problem has the following decision variables and associated dimensions, where $n_s$ is the number of client satellites. 
\begin{align*}
    \mathrm{Satellite~order~variables:}&~z_{i,j} \in \{0,1\}, &i,j \in [n_s], \\
    \mathrm{Orbit~radii:}&~r_{\mathrm{orbit}, i} \in [r_{\rm{orbit,min}}, r_{\rm{orbit,max}}], &i \in [n_s-1], \\ 
    \mathrm{Orbital~periods:}&~ T_{\rm{orbit}, i} \in [T_{\rm{orbit,min}}, T_{\rm{orbit,max}}], &i \in [n_s-1], \\
    \mathrm{Orbital~period~differences:}&~\Delta T_{\rm{orbit}, i} \in [\Delta T_{\rm{min}}, \Delta T_{\rm{max}}], &i \in [n_s-1], \\
    \mathrm{True~anomalies:}&~\theta_i \in [-\pi, \pi], &i \in [n_s-1], \\ 
    \mathrm{Transfer~times:}&~ t_{\rm{transfer}, i} \in [0, t_{\rm{transfer,max}}], &i \in [n_s-1], \\
    \mathrm{Maneuver~times:}&~ t_{\rm{maneuver}, i} \in \R^+, &i \in [n_s-1], \\
    \mathrm{Orbital~revolutions:}&~N_{\mathrm{orbit}, i} \in [50, 500], &i \in [n_s-1], \\ 
    \mathrm{Orbital~entry~mass~ratios:}&~f_{\rm{entry}, i} \in [1, 1.0025],  &i \in [n_s-1], \\ 
    \mathrm{Orbital~exit~mass~ratios:}&~f_{\rm{exit}, i} \in [1, 1.0025],  &i \in [n_s-1], \\ 
    \mathrm{Wet~mass:}&~m_{\rm{wet}} \in [m_{\rm{dry}}, 2000], &  \\ 
    \mathrm{Intermediate~masses:}&~m_{i,j} \in [m_{\rm{dry}}, 2000], &i \in [n_s-1],~j \in [5], \\
    \mathrm{Transferred~fuel~masses:}&~m_{\rm{fuel}, i} \in [m_{\rm{fuel,min}}, m_{\rm{fuel,max}}], &i \in [n_s]. 
\end{align*}

The objective function is to minimize the wet (i.e. fueled) mass of the satellite. Note that the orbital quantities define the phasing orbits that the servicer uses to transfer between client satellites, and all altitudes are converted to radii with respect to the center of the Earth for simplicity. 

The bounds $r_{\rm{orbit,min}},~r_{\rm{orbit,max}},~m_{\rm{dry}}$ are defined in Table~\ref{tab:oosparams} as the minimum and maximum servicer altitudes, and the servicer dry mass respectively. $m_{\rm{fuel,min}}$ and $m_{\rm{fuel,max}}$ are the minimum and maximum of the fuel requirements shown in Figure~\ref{fig:clientSatFuel}. Since $t_{\rm{maneuver}, i}$ is not in any nonlinear constraints, it doesn't require bounds. The remaining bounds are defined as a function of problem parameters such as specific impulse $I_{\rm{sp}}$, maximum service time $t_{\rm{max}}$, client orbital altitude $r_{\rm{client}}$ and client fuel requirements $\Delta m_{\rm{cf}, i},~i \in [n_s]$; as well as physical constants such as the gravitational constant $\mu$, and gravitational acceleration $g$.
\begin{align*}
    T_{\rm{client}} &= 2\pi\sqrt{\frac{r_{\rm{client}}}{\mu}} \\ 
    T_{\rm{orbit, min}} &= 2\pi\sqrt{\frac{r_{\rm{orbit, min}}}{\mu}} \\ 
    T_{\rm{orbit, max}} &= 2\pi\sqrt{\frac{r_{\rm{orbit, max}}}{\mu}} \\
    \Delta T_{\rm{orbit, min}} &= -\rm{max}(|T_{\rm{orbit},i} - T_{\rm{client}}|,~\forall i \in [n_s-1]) \\
    \Delta T_{\rm{orbit, max}} &= \rm{max}(|T_{\rm{orbit},i} - T_{\rm{client}}|,~\forall i \in [n_s-1]) \\
    t_{\rm{transfer,max}} &= 2\pi\sqrt{\frac{r_{\rm{orbit, max}} + r_{\rm{client}}}{8\mu}} \\ 
\end{align*}

\subsubsection{Linear Constraints.}

The constraints are given below with brief descriptions.

\resizebox{0.95\hsize}{!}{%
\begin{math}
\begin{aligned}
    \mathrm{Each~client~visited~once:} ~& \sum_{i=1}^{n_s} z_{i,j} = 1,~&\forall j \in [n_s] \\
    \mathrm{One~refuel~per~rendezvous:} ~&\sum_{j=1}^{n_s} z_{i,j} = 1,~&\forall i \in [n_s] \\
    \mathrm{Fuel~required~for}~i \rm{th~client:}~&m_{\rm{fuel}, i} = \sum_{j=1}^{n_s} \Delta m_{\rm{cf}, j} z_{i, j},~&\forall i \in [n_s]\\
    \mathrm{True~anomaly~from~client}~i~\mathrm{to}~i+1: ~&\theta_i = \sum_{j=1}^{n_s}\Bigg((-\pi + 2\pi j / n_s) (z_{i+1, j} - z_{i,j})\Bigg),~&\forall i \in [n_s-1]\\
    \mathrm{Wet~mass:}~& m_{\rm{wet}} = m_{1,1} + m_{\rm{fuel, 1}} \\
    \mathrm{Intermediate~fuel~transfers:}~& m_{i, 5} = m_{i+1, 1} + m_{\rm{fuel}, i+1},~&\forall i \in [n_s-2]\\
    \mathrm{Dry~mass:} ~& m_{n_s-1, 5} = m_{\rm{dry}} + m_{\rm{fuel}, n_s} \\
    \mathrm{Orbital~period~difference:} ~& \Delta T_{\rm{orbit}, i} = T_{\rm{orbit}, i} - T_{\rm{client}},~&\forall i \in [n_s-1]\\
    \mathrm{Total~maneuver~time:} ~& \sum_{i=1}^{n_s-1} t_{\rm{maneuver},i} \leq t_{\rm{max}} & \\ 
\end{aligned}
\end{math}}

\subsubsection{Nonlinear Constraints.}

The nonlinear constraints fall into 7 distinct forms, which are repeated in the satellite dynamical system. The list of 60 nonlinear constraints, as well as brief descriptions are below:

\begin{itemize}
    \item Transfer orbit entry burn ($n_s-1$ constraints): Describes mass ratio (entry mass over exit mass) of the satellite during transfer orbit entry. 
    \begin{equation*}
        \begin{split}
            f_{\rm{entry}, i} = \max\Bigg[&\exp\Bigg(\frac{1}{gI_{\rm{sp}}} \sqrt{\frac{\mu}{r_{\rm{orbit},i}}} \Bigg(\sqrt{\frac{2r_{\rm{client}}}{r_{\rm{client}} + r_{\rm{orbit}, i}}} - 1 \Bigg)\Bigg), \\
            &\exp\Bigg(\frac{1}{gI_{\rm{sp}}} \sqrt{\frac{\mu}{r_{\rm{client}}}} \Bigg(\sqrt{\frac{2r_{\rm{orbit},i}}{r_{\rm{client}} + r_{\rm{orbit}, i}}} - 1 \Bigg)\Bigg)\Bigg],~i \in [n_s-1].
        \end{split}
        \end{equation*}
    \item Transfer orbit exit burn ($n_s-1$ constraints): Describes the mass ratio (entry mass over exit mass) of the satellite during transfer orbit exit. 
    \begin{equation*}
    \begin{split}
        f_{\rm{exit}, i} = \max\Bigg[&\exp\Bigg(\frac{1}{gI_{\rm{sp}}} 
        \sqrt{\frac{\mu}{r_{\rm{client}}}} \Bigg(1 - \sqrt{\frac{2r_{\rm{orbit},i}}{r_{\rm{client}} + r_{\rm{orbit}, i}}}\Bigg)\Bigg), \\
        &\exp\Bigg(\frac{1}{gI_{\rm{sp}}} \sqrt{\frac{\mu}{r_{\rm{orbit},i}}} \Bigg(1 - \sqrt{\frac{2r_{\rm{client}}}{r_{\rm{client}} + r_{\rm{orbit}, i}}}\Bigg)\Bigg)\Bigg],~i \in [n_s-1].
    \end{split}
    \end{equation*}
    \item Mass conservation ($4(n_s-1)$ constraints): Couples the fractional change in mass of the satellite to the absolute change in mass during each burn phase. 
    \begin{align*}
        m_{i,1} &= f_{\rm{entry},i} m_{i, 2},~&i \in [n_s-1],\\
        m_{i,2} &= f_{\rm{exit},i} m_{i, 3},~&i \in [n_s-1],\\
        m_{i,3} &= f_{\rm{exit},i} m_{i, 4},~&i \in [n_s-1],\\
        m_{i,4} &= f_{\rm{entry},i} m_{i, 5},~&i \in [n_s-1].
    \end{align*}
    \item Phasing orbit period ($n_s-1$ constraints): Describes the period of the phasing orbit. 
    \begin{equation*}
        T_{\rm{orbit}, i} = 2 \pi \sqrt{\frac{r_{\rm{orbit},i}^3}{\mu}},~i \in [n_s-1].
    \end{equation*}
    \item Transfer time ($n_s-1$ constraints): Describes the Hohmann transfer time from the client to phasing orbit. 
    \begin{equation*}
        t_{\rm{transfer}, i} = 2 \pi \sqrt{\frac{(r_{\rm{client}} + r_{\rm{orbit},i})^3}{8 \mu}},~i \in [n_s-1].
    \end{equation*}
    \item Number of transfer orbit revolutions ($n_s-1$ constraints): Describes the number of revolutions in phasing orbit.  
    \begin{equation*}
        N_{\rm{orbit}, i} \Delta T_{\rm{orbit}, i} = T_{\rm{client}, i} \theta_i,~i \in [n_s-1].
    \end{equation*}
    \item Maneuver time ($n_s-1$ constraints): Describes the maneuver time (transfer and phasing time) between clients. 
    \begin{equation*}
        t_{\rm{maneuver},i} = t_{\rm{transfer}, i} + N_{\rm{orbit}, i} T_{\rm{orbit},i},~i \in [n_s-1].
    \end{equation*}
\end{itemize}
\end{APPENDICES}

% Acknowledgments here
\ACKNOWLEDGMENT{We thank Michael Luu for his help in the formulation of the satellite \gls{oos} problem in Section~\ref{sec:oos}, based on his previous work~[\cite{Luu2020}]. 

This research was supported in part by the Defense Advanced Research Projects Agency (DARPA), under contract HR0011-17-2-0028. The views, opinions and/or findings expressed are those of the author and should not be interpreted as representing the official views or policies of the Department of Defense or the U.S. Government.}

\bibliographystyle{informs2014}
\bibliography{main}

\begin{thebibliography}{37}
\providecommand{\natexlab}[1]{#1}
\providecommand{\url}[1]{\texttt{#1}}
\providecommand{\urlprefix}{URL }

\bibitem[{Bates et~al.(2003)Bates, Sienz, \protect\BIBand{}
  Langley}]{Bates2003}
Bates SJ, Sienz J, Langley DS (2003) {Formulation of the Audze-Eglais Uniform
  Latin Hypercube design of experiments}. \emph{Advances in Engineering
  Software} 34(8):493--506,
  \urlprefix\url{http://dx.doi.org/10.1016/S0965-9978(03)00042-5}.

\bibitem[{Bates et~al.(2004)Bates, Sienz, \protect\BIBand{}
  Toropov}]{Bates2004}
Bates SJ, Sienz J, Toropov VV (2004) {Formulation of the Optimal Latin
  Hypercube Design of Experiments Using a Permutation Genetic Algorithm}.
  \emph{Collection of Technical Papers - AIAA/ASME/ASCE/AHS/ASC Structures,
  Structural Dynamics and Materials Conference} 7:5217--5223,
  \urlprefix\url{http://dx.doi.org/10.2514/6.2004-2011}.

\bibitem[{Bergamini et~al.(2008)Bergamini, Grossmann, Scenna, \protect\BIBand{}
  Aguirre}]{Bergamini2008}
Bergamini ML, Grossmann I, Scenna N, Aguirre P (2008) {An improved piecewise
  outer-approximation algorithm for the global optimization of MINLP models
  involving concave and bilinear terms}. \emph{Computers and Chemical
  Engineering} 32(3):477--493,
  \urlprefix\url{http://dx.doi.org/10.1016/j.compchemeng.2007.03.011}.

\bibitem[{Bertsimas \protect\BIBand{} Dunn(2017)}]{Bertsimas2017}
Bertsimas D, Dunn J (2017) { Optimal classification trees }. \emph{Machine
  Learning} 106(7):1039--1082,
  \urlprefix\url{http://dx.doi.org/10.1007/s10994-017-5633-9}.

\bibitem[{Bertsimas \protect\BIBand{} Dunn(2019)}]{Bertsimas2019}
Bertsimas D, Dunn J (2019) \emph{{ Machine Learning Under a Modern Optimization
  Lens }} (Dynamic Ideas Press).

\bibitem[{Bertsimas \protect\BIBand{} Stellato(2021)}]{Bertsimas2021}
Bertsimas D, Stellato B (2021) {The voice of optimization}. \emph{Machine
  Learning} 110(2):249--277,
  \urlprefix\url{http://dx.doi.org/10.1007/s10994-020-05893-5}.

\bibitem[{Biggs et~al.(2017)Biggs, Hariss, \protect\BIBand{}
  Perakis}]{Biggs2017}
Biggs M, Hariss R, Perakis G (2017) {Optimizing Objective Functions Determined
  from Random Forests}. \emph{SSRN Electronic Journal} 1--46,
  \urlprefix\url{http://dx.doi.org/10.2139/ssrn.2986630}.

\bibitem[{Breiman et~al.(1984)Breiman, Friedman, Olshen, \protect\BIBand{}
  Stone}]{Breiman1984}
Breiman L, Friedman JH, Olshen RA, Stone CJ (1984) \emph{{Classification and
  Regression Trees}} (Taylor \& Francis), ISBN 9780412048418.

\bibitem[{Bussieck et~al.(2003)Bussieck, Drud, \protect\BIBand{}
  Meeraus}]{Bussieck2003}
Bussieck MR, Drud AS, Meeraus A (2003) {MINLPLib - A collection of test models
  for mixed-integer nonlinear programming}. \emph{INFORMS Journal on Computing}
  15(1):114--119,
  \urlprefix\url{http://dx.doi.org/10.1287/ijoc.15.1.114.15159}.

\bibitem[{Cortes et~al.(1995)Cortes, Jackel, \protect\BIBand{}
  Chiang}]{Cortes1995}
Cortes C, Jackel L, Chiang WP (1995) {Limits on Learning Machine Accuracy
  Imposed by Data Quality}. \emph{Proceedings of the Association for the
  Advancement of Artificial Intelligence Conference on Knowledge Discovery and
  Data Mining} 57--62.

\bibitem[{Drud(1994)}]{ArneS.Drud1994}
Drud AS (1994) {CONOPT - A Large Scale GRG Code}. \emph{ORSA Journal on
  Computing} 6(2).

\bibitem[{Dunning et~al.(2017)Dunning, Huchette, \protect\BIBand{}
  Lubin}]{DunningHuchetteLubin2017}
Dunning I, Huchette J, Lubin M (2017) Jump: A modeling language for
  mathematical optimization. \emph{SIAM Review} 59(2):295--320,
  \urlprefix\url{http://dx.doi.org/10.1137/15M1020575}.

\bibitem[{Duran \protect\BIBand{} Grossmann(1986)}]{Duran1986}
Duran MA, Grossmann IE (1986) {An Outer-Approximation Algorithm for a Class of
  Mixed-Integer Nonlinear Programs}. \emph{Mathematical Programming}
  36:307--339.

\bibitem[{Frazier(2018)}]{Frazier2018}
Frazier PI (2018) {Bayesian Optimization}. \emph{INFORMS TutORials in
  Operations Research} 9--11.

\bibitem[{Gambella et~al.(2021)Gambella, Ghaddar, \protect\BIBand{}
  Naoum-Sawaya}]{Gambella2021}
Gambella C, Ghaddar B, Naoum-Sawaya J (2021) {Optimization problems for machine
  learning: A survey}. \emph{European Journal of Operational Research}
  290(3):807--828,
  \urlprefix\url{http://dx.doi.org/10.1016/j.ejor.2020.08.045}.

\bibitem[{Gastegger et~al.(2017)Gastegger, Behler, \protect\BIBand{}
  Marquetand}]{Gastegger2017}
Gastegger M, Behler J, Marquetand P (2017) {Machine learning molecular dynamics
  for the simulation of infrared spectra}. \emph{Chemical Science}
  8(10):6924--6935, \urlprefix\url{http://dx.doi.org/10.1039/c7sc02267k}.

\bibitem[{Golinski(1970)}]{Golinski1970}
Golinski J (1970) {Optimal Synthesis Problems Solved by Means of Nonlinear
  Programming and Random Methods}. \emph{Journal of Mechanisms} 5(March
  1969):287--309.

\bibitem[{{Gurobi Optimization, LLC}(2021)}]{gurobi2021}
{Gurobi Optimization, LLC} (2021) {Gurobi Optimizer Reference Manual}.

\bibitem[{Horst et~al.(1989)Horst, Thoai, \protect\BIBand{} Tuy}]{Horst1989}
Horst R, Thoai N, Tuy H (1989) {On an Outer Approximation Concept in Global
  Optimization}. \emph{Optimization} 20(3):255--264,
  \urlprefix\url{http://dx.doi.org/10.1080/02331938908843440}.

\bibitem[{{Interpretable AI, LLC}(2022)}]{InterpretableAI}
{Interpretable AI, LLC} (2022) {Interpretable AI Documentation}.
  \urlprefix\url{https://www.interpretable.ai}.

\bibitem[{Kochkov et~al.(2021)Kochkov, Smith, Alieva, Wang, Brenner,
  \protect\BIBand{} Hoyer}]{Kochkov2021}
Kochkov D, Smith JA, Alieva A, Wang Q, Brenner MP, Hoyer S (2021) {Machine
  learning – accelerated computational fluid dynamics}. \emph{Proceedings of
  the National Academy of Sciences} 118,
  \urlprefix\url{http://dx.doi.org/10.1073/pnas.2101784118}.

\bibitem[{Lin \protect\BIBand{} Tsai(2012)}]{Lin2012}
Lin MH, Tsai JF (2012) {Range reduction techniques for improving computational
  efficiency in global optimization of signomial geometric programming
  problems}. \emph{European Journal of Operational Research} 216(1):17--25,
  \urlprefix\url{http://dx.doi.org/10.1016/j.ejor.2011.06.046}.

\bibitem[{Luu \protect\BIBand{} Hastings(2020)}]{Luu2020}
Luu M, Hastings D (2020) {Valuation of On-Orbit Servicing in Proliferated
  Low-Earth Orbit Constellations}. \emph{Proceedings of AIAA ASCEND 2020}
  0--14.

\bibitem[{Maragno et~al.(2021)Maragno, Wiberg, Bertsimas, Birbil, den Hertog,
  \protect\BIBand{} Fajemisin}]{Maragno2021}
Maragno D, Wiberg H, Bertsimas D, Birbil SI, den Hertog D, Fajemisin A (2021)
  {Mixed-Integer Optimization with Constraint Learning}. \emph{arXiv} 1--48.

\bibitem[{McKay et~al.(1979)McKay, Beckman, \protect\BIBand{}
  Conover}]{McKay1979}
McKay MD, Beckman RJ, Conover WJ (1979) {A comparison of three methods for
  selecting values of input variables in the analysis of output from a computer
  code}. \emph{Technometrics} 21(2):239--245,
  \urlprefix\url{http://dx.doi.org/10.1080/00401706.2000.10485979}.

\bibitem[{Mi{\v{s}}i{\'{c}}(2020)}]{Misic2020}
Mi{\v{s}}i{\'{c}} VV (2020) {Optimization of tree ensembles}. \emph{Operations
  Research} 68(5):1605--1624,
  \urlprefix\url{http://dx.doi.org/10.1287/opre.2019.1928}.

\bibitem[{Morawietz \protect\BIBand{} Artrith(2021)}]{Morawietz2021}
Morawietz T, Artrith N (2021) {Machine learning-accelerated quantum
  mechanics-based atomistic simulations for industrial applications}.
  \emph{Journal of Computer-Aided Molecular Design} 35(4):557--586,
  \urlprefix\url{http://dx.doi.org/10.1007/s10822-020-00346-6}.

\bibitem[{Ray(2003)}]{Ray2003}
Ray T (2003) {Golinski's speed reducer problem revisited}. \emph{AIAA Journal}
  41(3):556--558, \urlprefix\url{http://dx.doi.org/10.2514/2.1984}.

\bibitem[{Ryoo \protect\BIBand{} Sahinidis(1996)}]{Ryoo1996}
Ryoo HS, Sahinidis NV (1996) {A branch-and-reduce approach to global
  optimization}. \emph{Journal of Global Optimization} 8(2):107--138,
  \urlprefix\url{http://dx.doi.org/10.1007/bf00138689}.

\bibitem[{Sahinidis(1996)}]{Sahinidis1996}
Sahinidis NV (1996) {BARON: A general purpose global optimization software
  package}. \emph{Journal of Global Optimization} 8(2):201--205,
  \urlprefix\url{http://dx.doi.org/10.1007/bf00138693}.

\bibitem[{Sahinidis(2017)}]{sahinidis:baron:21.1.13}
Sahinidis NV (2017) \emph{{BARON 21.1.13: Global Optimization of Mixed-Integer
  Nonlinear Programs, {\em User's Manual}}}.

\bibitem[{Shewry \protect\BIBand{} Wynn(1987)}]{Shewry1987}
Shewry MC, Wynn HP (1987) {Maximum entropy sampling}. \emph{Journal of Applied
  Statistics} 14(2):165--170,
  \urlprefix\url{http://dx.doi.org/10.1080/02664768700000020}.

\bibitem[{Sun et~al.(2020)Sun, Cao, Zhu, \protect\BIBand{} Zhao}]{Sun2020}
Sun S, Cao Z, Zhu H, Zhao J (2020) {A Survey of Optimization Methods from a
  Machine Learning Perspective}. \emph{IEEE Transactions on Cybernetics}
  50(8):3668--3681,
  \urlprefix\url{http://dx.doi.org/10.1109/TCYB.2019.2950779}.

\bibitem[{Tagliarini et~al.(1991)Tagliarini, Christ, \protect\BIBand{} {Page,
  Edward}}]{Tagliarini1991}
Tagliarini GA, Christ JF, {Page, Edward} W (1991) {Optimization Using Neural
  Networks}. \emph{IEEE Transactions on Computers} 40(12):1347--1358.

\bibitem[{Verma(2000)}]{Verma2000}
Verma A (2000) {An introduction to automatic differentiation}. \emph{Current
  Science} 78(7):804--807,
  \urlprefix\url{http://dx.doi.org/10.1002/pamm.200310012}.

\bibitem[{Vielma(2015)}]{Vielma2015}
Vielma JP (2015) {Mixed Integer Linear Programming Formulation Techniques}.
  \emph{SIAM Review} 57(1):3--57,
  \urlprefix\url{http://dx.doi.org/10.1137/130915303}.

\bibitem[{W{\"{a}}chter \protect\BIBand{} Biegler(2006)}]{Wachter2006}
W{\"{a}}chter A, Biegler LT (2006) {On the implementation of an interior-point
  filter line-search algorithm for large-scale nonlinear programming}.
  \emph{Math. Program.} 106:25--57,
  \urlprefix\url{http://dx.doi.org/10.1007/s10107-004-0559-y}.

\end{thebibliography}

%% Here starts the e-companion (EC)
%%%%%%%%%%%%%%%%%%%%%%%%%%%%%%%%%%%%%%%%%%%%%%%%%%%%%%%%%%
% \ECSwitch

%\ECDisclaimer
%%%%%%%%%%%%%%%%%%%%%%%%%%%%%%%%%%%%%%%%%%%%%%%%%%%%%%%%%%

%%% Main head for the e-companion
% \ECHead{Proofs of Statements}

%%%%%%%%%%%%%%%%%
\end{document}